\documentclass[a4paper,fleqn]{cas-sc}

\usepackage[numbers]{natbib}

\def\tsc#1{\csdef{#1}{\textsc{\lowercase{#1}}\xspace}}
\tsc{WGM}
\tsc{QE}

\usepackage{amssymb}
\usepackage{latexsym}
\usepackage{amsmath}
\usepackage{subcaption}
\usepackage{float}
\usepackage{siunitx}

\usepackage{todonotes}

\DeclareMathOperator{\dive}{\nabla\cdot}

\begin{document}
\let\WriteBookmarks\relax
\def\floatpagepagefraction{1}
\def\textpagefraction{.001}

\shorttitle{IMEX-DG solver for atmosphere dynamics}    

\shortauthors{G. Orlando et al.}  

\title [mode = title]{An IMEX-DG solver for atmospheric dynamics simulations with adaptive mesh refinement}  

\author[1]{Giuseppe Orlando}
\ead{giuseppe.orlando@polimi.it}

\author[2]{Tommaso Benacchio}
\ead{tommaso.benacchio.ext@leonardo.com}

\author[1]{Luca Bonaventura}

\cormark[1]
\ead{luca.bonaventura@polimi.it}

\affiliation[1]{organization={Dipartimento di Matematica, Politecnico di Milano},
	addressline={Piazza Leonardo da Vinci 32}, 
	city={Milano},
	postcode={20133}, 
	country={Italy}}

\affiliation[2]{organization={Future Rotorcraft Technologies, Leonardo S.p.a.},
	addressline={Cascina Costa}, 
	city={Samarate (VA)},
	postcode={21017}, 
	country={Italy}}

\cortext[1]{Corresponding author}

\begin{abstract}
We present an accurate and efficient solver for atmospheric dynamics simulations that allows for non-conforming mesh refinement. The model equations are the conservative Euler equations for compressible flows. The numerical method is based on an $h-$adaptive Discontinuous Galerkin spatial discretization and on a second order Additive Runge Kutta IMEX method for time discretization, especially designed for low Mach regimes. The solver is implemented in the framework of the \textit{deal.II} library, whose mesh refinement capabilities are employed to enhance efficiency. A number of numerical experiments based on classical benchmarks for atmosphere dynamics demonstrate the properties and advantages of the proposed method.	
\end{abstract}




\maketitle

\section{Introduction}
\label{sec:intro} \indent

The fully compressible equations of gas dynamics have become the standard of mathematical models of the atmosphere in the last couple of decades, see e.g. the literature review in \cite{steppeler:2003}. The efficient numerical solution of these equations in presence of gravitational forcing is crucial for weather and climate predictions and poses several computational challenges.

The Discontinuous Galerkin (DG) discretization method has been successfully applied to atmospheric modelling, see e.g. \cite{giraldo:2010, kelly:2012} or the review in \cite{bonaventura:2012}. Discontinuous finite element $h-$adaptive approaches were proposed in \cite{kopera:2014, muller:2013, yelash:2014}, while $p-$adaptive DG methods for Numerical Weather Prediction (NWP) were introduced in \cite{tumolo:2015}. However, severe time step stability restrictions may be required if the DG method is coupled with explicit time discretization schemes. The use of implicit and semi-implicit methods allows to relax these restrictions and employ much longer time steps, as recognized since a long time in the literature on atmosphere dynamics and numerical weather prediction, see among many others \cite{benacchio:2014, bonaventura:2000, giraldo:2005, giraldo:2013, melvin:2019, smolarkiewicz:2019}.

In this work, we will present an application of the $h-$adaptive DG method recently proposed in \cite{orlando:2022b} to classical benchmarks for atmospheric flow modelling. The proposed method is based on an accurate and flexible DG space discretization and an IMplicit-EXplicit (IMEX) time discretization, whose properties and theoretical justifications are discussed in detail in \cite{orlando:2022b} and only briefly summarized in the following. The adaptive discretization is implemented in the framework of the numerical library \textit{deal.II} \cite{arndt:2022, bangerth:2007}, an open-source environment, which provides non-conforming $h-$refinement capabilities that are exploited here for the first time in the numerical simulation of atmospheric flows. 
The development of a specific NWP code is thus embedded into a framework continuously adapted to novel architectures and extended to include novel versions of the DG method. As a result, the portability and the continuous development of the code are greatly enhanced. The results achieved shows that a dynamical core for high resolution, fine scale atmospheric modelling can be based on a widely used open source finite element library. It will also be shown that the proposed adaptive implementation can achieve good parallel efficiency, even without employing specific optimizations. 

The paper is structured as follows. The model equations are reviewed in Section \ref{sec:modeleq}. The time and space discretizations are briefly outlined and discussed in Section \ref{sec:disc}. The validation of the proposed method and its applications to a number of significant benchmarks with and without orography are reported in Sections \ref{sec:tests_without}, \ref{sec:tests_with}, respectively. Some conclusions and considerations about open issues and future work are presented in Section \ref{sec:conclu}.

\section{The model equation}
\label{sec:modeleq} \indent

Let \(\Omega \subset \mathbb{R}^{d}, 2 \le d \le 3\) be a connected open bounded set with a sufficiently smooth boundary \(\partial\Omega\) and denote by \(\mathbf{x}\) the spatial coordinates and by \(t\) the temporal coordinate. We consider the classical unsteady compressible Euler equations, written in flux form as:
\begin{eqnarray}
\label{eq:euler_comp}
\frac{\partial \rho}{\partial t}  
+ \dive\left(\rho\mathbf{u}\right) &=& 0 \nonumber \\
\frac{\partial\left(\rho\mathbf{u}\right)}{\partial t}  
+ \dive\left(\rho\mathbf{u} \otimes \mathbf{u}\right) 
+ \nabla p &=& -\rho g\mathbf{k}  \\
\frac{\partial (\rho E)}{\partial t}  
+ \dive\left[\left(\rho E + p\right)\mathbf{u} \right] &=&  -\rho g\mathbf{k} \cdot \mathbf{u} \nonumber
\end{eqnarray}
for $ \mathbf{x} \in \Omega, $ $t \in [0, T_f], $ supplied with suitable initial and boundary conditions. Here $ T_f $ is the final time, $ \rho $ is the density, $ \mathbf{u} $ is the fluid velocity, $ p $ is the pressure,  $ g = \SI{9.81}{\meter\per\second\squared}$ is the acceleration of gravity and $ \mathbf{k} $ is the upward pointing unit vector in the standard Cartesian reference frame. The total energy $ \rho E $ can be rewritten as $ \rho E = \rho e + \rho k,$ where $ e $ is the internal energy and $ k = \|\mathbf{u}\|^2/2 $ is the kinetic energy. We also introduce the specific enthalpy $ h = e + p/\rho $ and we notice that one can rewrite the energy flux as
$
\left(\rho E + p\right)\mathbf{u} = \left(e + k + p/\rho\right)\rho\mathbf{u} = \left(h+k\right)\rho\mathbf{u}.
$
The equations can then be rewritten as 
\begin{eqnarray}
\label{eq:euler_comp2}
\frac{\partial \rho}{\partial t} + \dive\left(\rho\mathbf{u}\right) &=& 0 \nonumber \\
\frac{\partial \left(\rho\mathbf{u}\right)}{\partial t} +
\dive\left(\rho \mathbf{u} \otimes \mathbf{u}\right) + 
\nabla p &=& -\rho g \mathbf{k} \\
\frac{\partial (\rho E)}{\partial t} + \dive\left[\left(h + k\right)\rho\mathbf{u} \right] &=& -\rho g \mathbf{k} \cdot\mathbf{u}.  \nonumber 
\end{eqnarray}
The above equations are complemented by the equation of state (EOS) for ideal gases, given by $ p = \rho RT, $ with \(R\) being the specific gas constant and $ T $  the absolute temperature.  As explained in \cite{orlando:2022b}, the proposed method allows to handle in the same framework also more general equations of state for real gases, thus opening the way for example to the inclusion of effects due to water vapour and moist species. In \cite{orlando:2022b} it is also shown how to extend the proposed approach to the viscous case by Strang splitting \cite{strang:1968}. This extension will be required by some of the numerical simulations described in Section \ref{sec:tests_without}. The potential temperature \(\theta\) is often preferred in atmospheric applications as the unknown in the energy equation. For an ideal gas, potential temperature is defined as
$ \theta = T \left(p_0/p\right)^{\frac{\gamma - 1}{\gamma}}, $
where \(p_0\) is a reference pressure and \(\gamma\) denotes the specific heats ratio. In the present work, we consider \(p_{0} = 10^{5} \hspace{0.05cm}\SI{}{\pascal}\). 
Notice that, even though the general existence, uniqueness and regularity of the solutions of the Euler equations \eqref{eq:euler_comp}, \eqref{eq:euler_comp2} is still an advanced research topic (see e.g. the results and references in \cite{feireisl:2017}) it is customary to give for granted existence, uniqueness and sufficient regularity of the solutions when dealing with atmospheric flows at low Mach numbers, such as those considered in this paper.

\section{The numerical framework}
\label{sec:disc} \indent

In the low Mach number limit, pressure gradient terms yield stiff components for the resulting semidiscretized ODE system, since the pressure gradients of the non dimensional equations are proportional to $1/M^2$, wuth $M$ denoting the Mach number, see e.g. the analysis in \cite{casulli:1984, munz:2003}. Therefore, following \cite{casulli:1984, dumbser:2016}, the method proposed in \cite{orlando:2022b} couples implicitly the energy equation to the momentum one, while the continuity equation is treated in a fully explicit fashion. For the time discretization, an IMplicit EXplicit (IMEX) Additive Runge Kutta method (ARK) \cite{kennedy:2003} method will be used. These methods are widely employed for ODE systems that include both stiff and non-stiff components, to which the implicit and explicit schemes are applied, respectively. Here, we consider a variant of the IMEX method proposed in \cite{giraldo:2013}, whose coefficients are presented in the Butcher tableaux reported in Table \ref{tab:ark2_butch} for the explicit and implicit method, respectively. We consider \(\chi = 2 - \sqrt{2}\), so that the implicit part of the IMEX scheme coincides with the TR-BDF2 method considered in \cite{tumolo:2015}. A solver based on this method for the incompressible Navier Stokes equations with an artificial compressibility formulation has been proposed in \cite{orlando:2022a}, highlighting the robustness of the proposed approach in the low Mach number limit. Notice also that, as discussed in \cite{orlando:2022b}, we take \(\alpha = 0.5\), rather than the value originally chosen in \cite{giraldo:2013}, in order to improve the monotonicity of the explicit part of the method.

\begin{table}[h!]
	\begin{center}
		\begin{tabular}{c|ccc}
			0 & 0 & &  \\
			$ \chi $ & $ \chi $ & 0 & \\
			1 & $ 1 - \alpha $ & $ \alpha $ & 0 \\
			\hline
			& $ \frac{1}{2} - \frac{\chi}{4} $ & $ \frac{1}{2} - \frac{\chi}{4} $ & $ \frac{\chi}{2} $
		\end{tabular}
		\hskip 2cm
		\begin{tabular}{c|ccc}
			0 & 0 & & \\
			$ \chi $ & $ \frac{\chi}{2}$ & $ \frac{\chi}{2} $ & \\
			1 & $ \frac{1}{2\sqrt{2}} $ & $ \frac{1}{2\sqrt{2}}$ & $ 1 - \frac{1}{\sqrt{2}} $ \\
			\hline
			& $ \frac{1}{2} - \frac{\chi}{4} $ & $ \frac{1}{2} - \frac{\chi}{4} $ & $ \frac{\chi}{2} $
		\end{tabular}
	\end{center}
	\caption{\it Butcher tableaux of the components of the ARK2 method: explicit method (left), implicit method (right)}
	\label{tab:ark2_butch}
\end{table}

We now briefly describe the application of this IMEX method outlined above to the time semi-discretization of equations \eqref{eq:euler_comp2}, the reader is referred to \cite{orlando:2022b} for the complete description of the scheme and for a complete convergence test on a benchmark with analytic solution. As it can be seen from the Butcher tableaux, the first stage is only formal, yielding $\rho^{(n,1)} =\rho^{n}$ for density and all the other prognostic quantities. 
The second stage is then given by
\begin{eqnarray}
\label{eq:euler_comp_hyp_tdisc1}
\rho^{(n,2)} &=& \rho^{n} - a_{21} \Delta t \dive\left(\rho^{n} \mathbf{u}^{n} \right) \nonumber \\
\rho^{(n,2)}\mathbf{u}^{(n,2)} &+& 
\tilde{a}_{22} \Delta t \nabla p^{(n,2)} = \mathbf{m}^{(n,2)} \\
\rho^{(n,2)} E^{(n,2)} &+& 
\tilde{a}_{22} \Delta t \dive\left(h^{(n,2)} \rho^{(n,2)} \mathbf{u}^{(n,2)} \right) = \mathbf{e}^{(n,2)},
\nonumber
\end{eqnarray}
where we have set
\begin{eqnarray}
\label{eq:euler_comp_hyp_tdisc1_rhss}
\mathbf{m}^{(n,2)} &=& \rho^{n} \mathbf{u}^{n}
- a_{21} \Delta t \dive\left(\rho^{n} \mathbf{u}^{n} \otimes\mathbf{u}^{n}\right) 
- \tilde{a}_{21} \Delta t\nabla p^{n}
- \tilde{a}_{21} \Delta t \rho^{n} g \mathbf{k} - \tilde a_{22} \Delta t \rho^{(n,2)} g \mathbf{k} \\
& \hskip 0.4cm & \nonumber \\
\mathbf{e}^{(n,2)} &=&\rho^{n} E^{n}  -\tilde a_{21} \Delta t \dive \left(h^{n} \rho^{n} \mathbf{u}^{n} \right)
- a_{21} \Delta t \dive\left(k^{n} \rho^{n} \mathbf{u}^{n} \right) 
- \tilde{a}_{21} \rho^{n} g \mathbf{k} \cdot \mathbf{u}^{n} - \tilde a_{22} \rho^{(n,2)} g\mathbf{k} \cdot \mathbf{u}^{(n,2)}. \nonumber 
\end{eqnarray}
For the third stage, one can write formally  
\begin{eqnarray}
\label{eq:euler_comp_hyp_tdisc2}
\rho^{(n,3)} &=& \rho^{n} - a_{31} \Delta t \dive\left(\rho^{n} \mathbf{u}^{n}\right)
- a_{32} \Delta t \dive\left(\rho^{(n,2)} \mathbf{u}^{(n,2)}\right)  \nonumber \\
\rho^{(n,3)}\mathbf{u}^{(n,3)} &+&
\tilde{a}_{33} \Delta t \nabla p^{(n,3)} = \mathbf{m}^{(n,3)}  \\
\rho^{(n,3)}E^{(n,3)} &+&
\tilde{a}_{33} \Delta t \dive\left(h^{(n,3)} \rho^{(n,3)} \mathbf{u}^{(n,3)}\right)= \mathbf{e}^{(n,3)}, \nonumber
\end{eqnarray}
where the right hand sides are defined as  
\begin{eqnarray}
\label{eq:euler_comp_tdisc2_rhss}
\mathbf{m}^{(n,3)} &=& \rho^{n} \mathbf{u}^n 
- a_{31} \Delta t \dive\left(\rho^{n} \mathbf{u}^{n} \otimes \mathbf{u}^{n}\right) 
- \tilde a_{31} \Delta t \nabla p^{n} \nonumber \\
&-& a_{32} \Delta t \dive\left(\rho^{(n,2)} \mathbf{u}^{(n,2)} \otimes \mathbf{u}^{(n,2)}\right) 
- \tilde a_{32} \Delta t \nabla p^{(n,2)} 
- \tilde{a}_{31} \Delta t g\rho^{n} \mathbf{k} - \tilde{a}_{32} \Delta t g\rho^{(n,2)} \mathbf{k} - \tilde{a}_{33} \Delta t g\rho^{(n,3)} \mathbf{k} \\
& \hskip 0.4cm & \nonumber \\
\mathbf{e}^{(n,3)} &=& \rho^n E^n -\tilde a_{31} \Delta t \dive \left(h^n  \rho^n\mathbf{u}^n \right)
- a_{31} \Delta t \dive\left(k^n \rho^n\mathbf{u}^n \right) 
- \tilde{a}_{31} \Delta t \rho^n g\mathbf{k} \cdot \mathbf{u}^n \nonumber \\
&-& \tilde a_{32} \Delta t \dive \left(h^{(n,2)} \rho^{(n,2)}\mathbf{u}^{(n,2)}\right)
- a_{32} \Delta t \dive\left(k^{(n,2)} \rho^{(n,2)}\mathbf{u}^{(n,2)}\right) 
- \tilde{a}_{32} \Delta t \rho^{(n,2)} g\mathbf{k} \cdot \mathbf{u}^{(n,2)} - \tilde{a}_{33} \Delta t \rho^{(n,3)} g\mathbf{k} \cdot \mathbf{u}^{(n,3)}. \nonumber 
\end{eqnarray}
In each stage, after spatial discretization has been performed, we formally substitute the momentum into the energy equation and we obtain an equation for the pressure \cite{orlando:2022b}. Notice that, in the above equations, the gravity source term is treated implicitly, using for this term exactly the same time discretization approach as for the pressure gradient terms, in order to avoid that the time discretizations induces spurious disturbances of hydrostatic balance.

For the spatial discretization, we consider a decomposition of the domain \(\Omega\) into a family of hexahedra \(\mathcal{T}_h\) (quadrilaterals in the two-dimensional case) and denote each element by \(K\). The skeleton \(\mathcal{E}\) denotes the set of all element faces and \(\mathcal{E} = \mathcal{E}^{I} \cup \mathcal{E}^{B}\), where \(\mathcal{E}^{I}\) and \(\mathcal{E}^{B}\) are the subset of interior and boundary faces, respectively. Suitable jump and average operators can then be defined as customary for finite element discretizations. A face \(\Gamma \in \mathcal{E}^{I}\) shares two elements that we denote by \(K^{+}\) with outward unit normal \(\mathbf{n}^{+}\) and \(K^{-}\) with outward unit normal \(\mathbf{n}^{-}\), whereas for a face \(\Gamma \in \mathcal{E}^{B}\) we simply denote by \(\mathbf{n}\) the outward unit normal.
For a scalar function \(\varphi\) the jump is defined as
$$\left[\left[\varphi\right]\right] = \varphi^{+}\mathbf{n}^{+} + \varphi^{-}\mathbf{n}^{-} \quad \text{if }\Gamma \in \mathcal{E}^{I} \qquad \left[\left[\varphi\right]\right] = \varphi\mathbf{n} \quad \text{if }\Gamma \in \mathcal{E}^{B}.$$
The average is defined as
$$\left\{\left\{\varphi\right\}\right\} = \frac{1}{2}\left(\varphi^{+} + \varphi^{-}\right) \quad \text{if }\Gamma \in \mathcal{E}^{I} \qquad \left\{\left\{\varphi\right\}\right\} = \varphi \quad \text{if }\Gamma \in \mathcal{E}^{B}.$$
Similar definitions apply for a vector function \(\boldsymbol{\varphi}\):
\begin{align*}
&\left[\left[\boldsymbol{\varphi}\right]\right] = \boldsymbol{\varphi}^{+}\cdot\mathbf{n}^{+} 
+\boldsymbol{\varphi}^{-}\cdot\mathbf{n}^{-} \quad \text{if }\Gamma \in \mathcal{E}^{I} \qquad 
\left[\left[\boldsymbol{\varphi}\right]\right] = \boldsymbol{\varphi}\cdot\mathbf{n} \quad \text{if }\Gamma \in \mathcal{E}^{B} \\
&\left\{\left\{\boldsymbol{\varphi}\right\}\right\} = \frac{1}{2}\left(\boldsymbol{\varphi}^{+} + \boldsymbol{\varphi}^{-}\right) \quad \text{if }\Gamma \in \mathcal{E}^{I} \qquad \left\{\left\{\boldsymbol{\varphi}\right\}\right\} = \boldsymbol{\varphi} \quad \text{if }\Gamma \in \mathcal{E}^{B}.
\end{align*}
For vector functions, it is also useful to define a tensor jump as:
$$\left<\left<\boldsymbol{\varphi}\right>\right> = \boldsymbol{\varphi}^{+}\otimes\mathbf{n}^{+} 
+ \boldsymbol{\varphi}^{-}\otimes\mathbf{n}^{-} \quad \text{if }\Gamma \in \mathcal{E}^{I} 
\qquad \left<\left<\boldsymbol{\varphi}\right>\right> = \boldsymbol{\varphi}\otimes\mathbf{n} \quad \text{if }\Gamma \in \mathcal{E}^{B}.$$
We introduce the following finite element spaces
\[Q_{r} = \left\{v \in L^2(\Omega) : v\rvert_K \in \mathbb{Q}_{r} \quad \forall K \in \mathcal{T}_h\right\} \ \ \
\mathbf{V}_{r} = \left[Q_{r}\right]^d,\]
where \(\mathbb{Q}_{r}\) is the space of polynomials of degree \(r\) in each coordinate direction and we recall that $d$ denotes the number of spatial dimensions. We then denote by \(\boldsymbol{\varphi}_i(\mathbf{x})\) the basis functions for the space \(\mathbf{V}_{r}\) and by \(\psi_i(\mathbf{x})\) the basis functions for the space \(Q_{r}\), the finite element spaces chosen for the discretization of the velocity and of the pressure (as well as the density), respectively, so that
\begin{equation*}
\mathbf{u}\approx \sum_{j = 1}^{\text{dim}(\mathbf{V}_{r})}u_j(t)\boldsymbol{\varphi}_j(\mathbf{x}) \qquad p \approx \sum_{j = 1}^{\text{dim}(Q_{r})}p_j(t)\psi_j(\mathbf{x}).
\end{equation*}
The shape functions are chosen to be the Lagrange interpolation polynomials for
the support points of \((r + 1)\)-order Gauss-Lobatto quadrature rule along each coordinate direction. Notice that the discretization does not change in principle for non-conforming grids, for which the features of the \textit{deal.II} library are exploited in our implementation. The use of tensor product finite element spaces allows to reduce all the issues to those concerning one-dimensional shape functions at one-dimensional quadrature points, while the extension is handled by means of tensor products. For face integrals at hanging nodes, the coarser of the two adjacent cells must interpolate the values to a subface (with evaluation points scaled either to the first or the second half of the reference interval representing the edge. The only constraint is that the refinement ratio between adjacent elements should be at most 1/2.

Given these definitions, the weak formulation for the momentum equation at each stage \(s = 2,3\) of the IMEX scheme can be written in compact form as
\begin{equation}
\mathbf{A}^{(n,s)}\mathbf{U}^{(n,s)} + \mathbf{B}^{(n,s)}\mathbf{P}^{(n,s)} = \mathbf{F}^{(n,s)},
\end{equation}
where we have set
\begin{eqnarray}
A_{ij}^{(n,s)} &=& \sum_K\int_K\rho^{(n,s)}\boldsymbol{\varphi}_j\cdot\boldsymbol{\varphi}_id\Omega \\
B_{ij}^{(n,s)} &=& \sum_K \int_K -\tilde{a}_{ss}\Delta t\nabla \cdot \boldsymbol{\varphi}_i \psi_j d\Omega 
+ \sum_{\Gamma\in\mathcal{E}}\int_{\Gamma} \tilde{a}_{ss} \Delta t\left\{\left\{\psi_j\right\}\right\}\left[\left[\boldsymbol{\varphi}_i\right]\right]d\Sigma
\end{eqnarray}
with \(\mathbf{U}^{(n,s)}\) denoting the vector of the degrees of freedom associated to the velocity field and \(\mathbf{P}^{(n,s)}\) denoting the vector of the degrees of freedom associated to the pressure. Moreover, \(\rho^{(n,s)}\) denotes the approximation of the density at stage \(s\) and \(\tilde{a}_{ss}\) is the corresponding coefficient of the Butcher tableaux in Table \ref{tab:ark2_butch} for the implicit part. Notice that the definitions introduced above entail that centered fluxes are used in the definition of the discrete pressure gradient, while an upwind flux has been employed in the definition of \(\mathbf{F}^{(n,s)}\) for the discretization of the terms computed at time level $n.$ Analogously, the weak formulation for the energy equation can be expressed as
\begin{equation}
\mathbf{C}^{(n,s)}\mathbf{U}^{(n,s)} = \mathbf{G}^{(n,s)},
\end{equation}
where we have set
\begin{equation}
C_{ij}^{(n,s)} = \sum_K \int_K -\tilde{a}_{ss} \Delta t h^{(n,s)}\rho^{(n,s)} \boldsymbol{\varphi}_j \cdot \nabla\psi_i d\Omega \nonumber \\ 
+ \sum_{\Gamma\in\mathcal{E}} \int_{\Gamma} \tilde{a}_{ss} \Delta t \left\{\left\{h^{(n,s)}\rho^{(n,s)} \boldsymbol{\varphi}_j\right\}\right\} \cdot \left[\left[\psi_i\right]\right] d\Sigma, 
\end{equation}
where \(h^{(n,s)}\) represents the approximation of the enthalpy at stage \(s\). Formally, we can then derive \(\mathbf{U}^{(n,s)} = (\mathbf{A}^{(n,s)})^{-1}\left(\mathbf{F}^{(n,s)} - \mathbf{B}^{(n,s)}\mathbf{P}^{(n,s)}\right)\) and obtain the following relation 
\begin{equation}
\mathbf{C}^{(n,s)}(\mathbf{A}^{(n,s)})^{-1} \left(\mathbf{F}^{(n,s)} - \mathbf{B}^{(n,s)}\mathbf{P}^{(n,s)}\right) = \mathbf{G}^{(n,s)}.
\end{equation}
Again, in the definition of \(\mathbf{G}^{(n,s)}\) an upwind flux has been employed for the discretization of the terms computed at time level $n.$  Moreover, taking into account that
$
\rho^{(n,s)}E^{(n,s)} = \rho^{(n,s)}e^{(n,s)}(p^{(n,s)}) + \rho^{(n,s)}k^{(n,s)},
$ 
we decompose \(\mathbf{G}^{(n,s)}\) and we finally obtain 
\begin{equation}
\mathbf{C}^{(n,s)}(\mathbf{A}^{(n,2)})^{-1}\left(\mathbf{F}^{(n,s)} - \mathbf{B}^{(n,s)}\mathbf{P}^{(n,s)}\right) = -\mathbf{D}^{(n,s)}\mathbf{P}^{(n,2)} +  \tilde{\mathbf{G}}^{(n,s)},
\end{equation}
where we have set
\begin{align}
D_{ij}^{(n,s)} = &\sum_K\int_K\rho^{(n,s)}e^{(n,s)}(\psi_j)\psi_id\Omega.
\end{align}
The above system can be solved in terms of \(\mathbf{P}^{(n,s)}\) according to the fixed point procedure described in \cite{dumbser:2016}. More specifically, setting
$\mathbf{P}^{(n,s,0)} = \mathbf{P}^{(n,s - 1)} $ and  $ k^{(n,s,0)} = k^{(n,s - 1)} $  
for \(l=1,\dots, M\) one solves the equation
\begin{equation}
\left(\mathbf{D}^{(n,s,l)} - \mathbf{C}^{(n,s,l)}(\mathbf{A}^{(n,s)})^{-1}
\mathbf{B}^{(n,s)}\right)\mathbf{P}^{(n,s,l+1)} =  
\tilde{\mathbf{G}}^{(n,s,l)} - \mathbf{C}^{(n,s,l)}(\mathbf{A}^{(n,s)})^{-1}\mathbf{F}^{(n,s,l)} \nonumber
\end{equation}
and updates the velocity solving
\begin{equation}
\mathbf{A}^{(n,s)}\mathbf{U}^{(n,s,l+1)} = \mathbf{F}^{(n,s,l)} - \mathbf{B}^{(n,s)}\mathbf{P}^{(n,s,l+1)}. \nonumber
\end{equation}
Once the iterations have been completed, one sets $ \mathbf{u}^{(n,s)} = \mathbf{u}^{(n,s,M+1)} $  and $ E^{(n,s)} $ accordingly. It is important to point out that the scheme outlined above only requires the solution of linear systems of a size equal to that of the number of discrete degrees of freedom associated to a scalar variable, as in \cite{dumbser:2016}, which is crucial for the overall efficiency. 

\section{Numerical results without orography}
\label{sec:tests_without} \indent

The numerical method outlined in the previous Sections has been validated in a number of benchmarks relevant for atmospheric applications. We define 
$${\cal H} = \min\{\mathrm{diam}(K) | K \in {\cal T}_{h} \} $$ 
and we define two Courant numbers, one based on the speed of sound denoted by \(C\), the so-called acoustic Courant number, and one based on the local velocity of the flow, the so-called advective Courant number, denoted by \(C_u\): 
\begin{equation}
\label{eq:Courant}
C = rc \Delta t/{\cal H}, \ \ \ \ \    C_u = rU \Delta t/{\cal H}
\end{equation}
where \(c\) is the magnitude of the speed of sound and \(U\) is the magnitude of the flow velocity. Notice that both the Courant numbers depend on the polynomial degree \(r\). In general, we consider \(\gamma = 1.4\) and \(R = \SI{287}{\joule\per\kilogram\per\kelvin}\) for all the simulations. In this Section, we will focus on the results obtained in benchmarks without orography, while we will discuss benchmarks with non trivial orographic profiles in the following Section \ref{sec:tests_with}. Notice also that, while Sections \ref{ssec:igw}, \ref{ssec:hydrostatic} and \ref{ssec:3D_mountain} include non adaptive results, which provide a general validation of the proposed discretization approach, Sections \ref{ssec:straka}, \ref{ssec:cold_bubble}, \ref{ssec:3D_bubble} and \ref{ssec:nonhydrostatic} are devoted to applications of the proposed adaptation technique.

\subsection{2D inertia-gravity waves}
\label{ssec:igw}

\begin{figure}[h!]
	\includegraphics[width=0.45\textwidth]{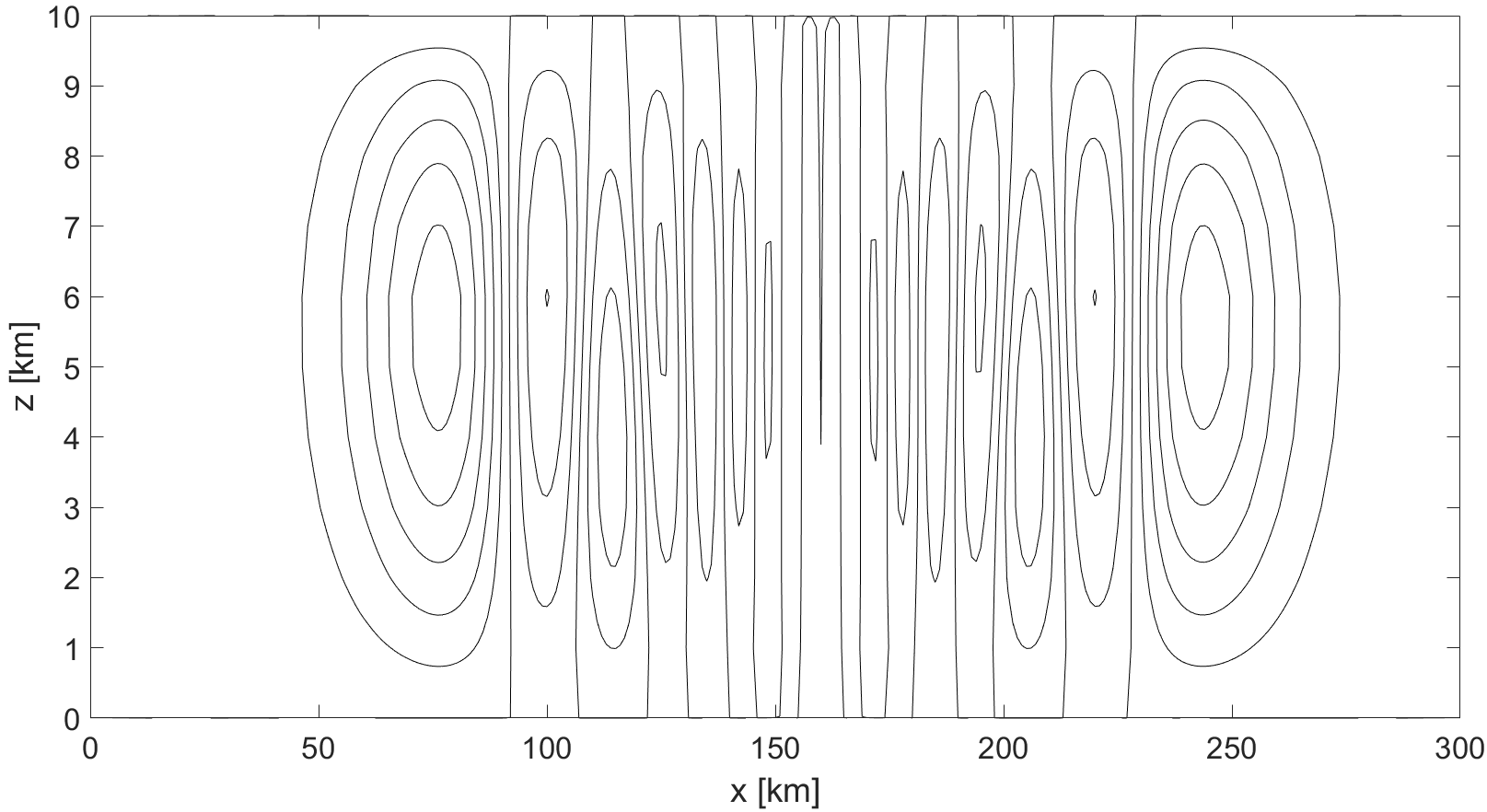} a)
	\includegraphics[width=0.45\textwidth]{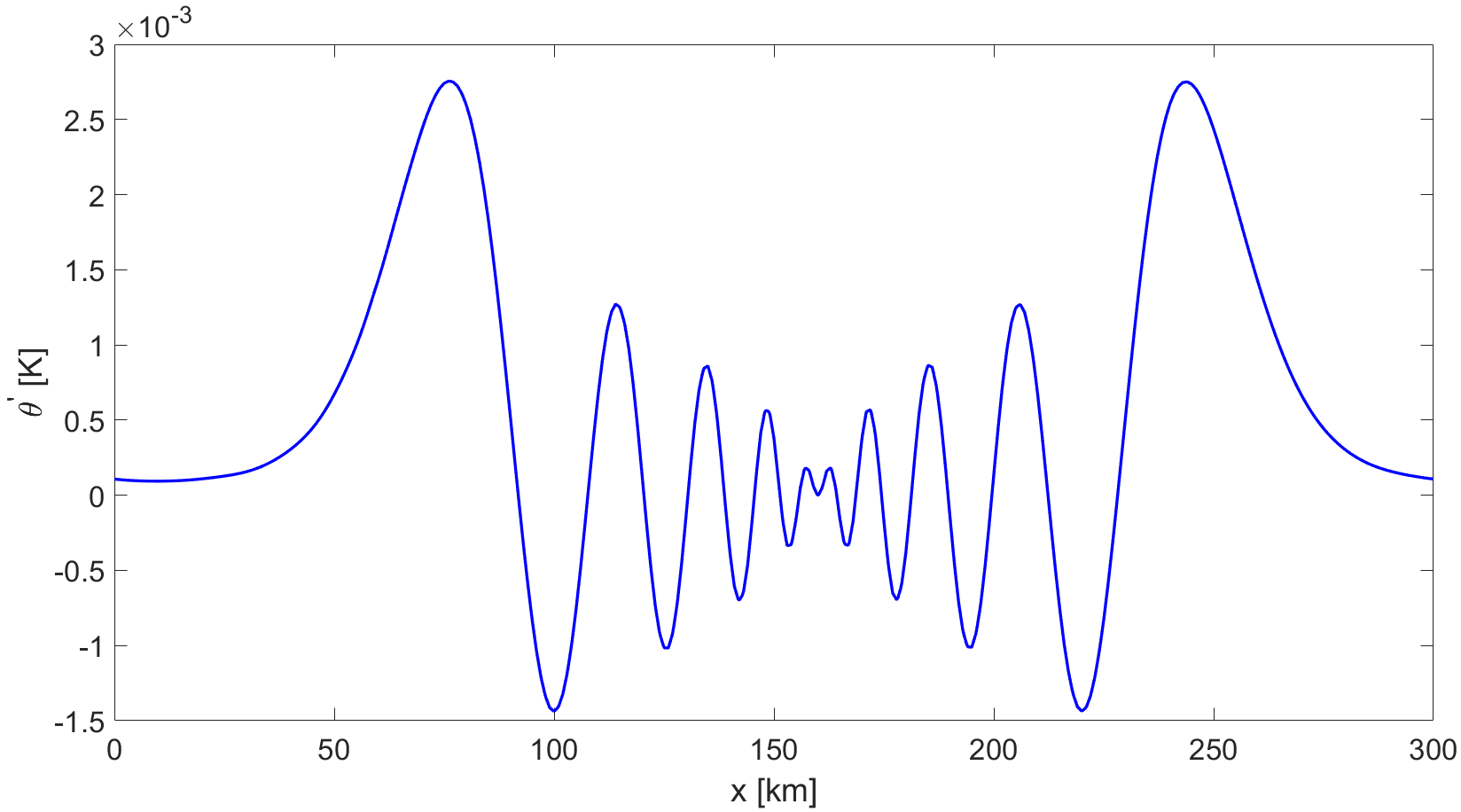} b)
	\caption{Inertia-gravity waves benchmark: a) potential temperature deviation  at \(t = \SI{3000}{\second}\), contours are plotted  from \(\SI{-0.0015}{\kelvin}\) to \(\SI{0.003}{\kelvin}\) with interval equal to \(\SI{5e-4}{\kelvin}\); b) potential temperature deviation at \(t = \SI{3000}{\second}\) along \(z = \SI{5}{\kilo\meter}\) height.}
	\label{fig:IGW}
\end{figure}

Inertia-gravity waves in a two-dimensional vertical section of the atmosphere constitute a classical benchmark for atmospheric flow models, see e.g. \cite{bonaventura:2000, melvin:2019, skamarock:1994}. In particular, we set the background potential temperature $\bar{\theta} = T_{ref}\exp\left({N^2 z}/{g}\right), $ where \(N = \SI{0.01}{\per\second}\) denotes the buoyancy frequency and \(T_{ref} = \SI{300}{\kelvin}\). The background density and pressure are defined as
$$ \bar{p} = \exp\left\{ 1 - \frac{g^2}{N^2}\frac{\gamma - 1}{\gamma}\frac{\rho_{ref}}{p_{ref}} \left[1 - \exp\left(-\frac{N^2 z}{g}\right)\right] \right\} \ \ \ \ \ \
\bar{\rho} = \rho_{ref}\left( \frac{p}{p_{ref}}\right)^{\frac{1}{\gamma}} \exp\left(-\frac{N^2 z}{g}\right) $$
with \(p_{ref} = 10^{5} \hspace{0.05cm} \SI{}{\pascal}\) and \(\rho_{ref} =  {p_{ref}}/{R T_{ref}}\). The domain is \(\Omega = \left(0, 300\right) \times \left(0, 10\right) \SI{}{\kilo\meter}\) and we consider the following perturbation for the potential temperature 
\begin{equation}
\theta^{'} = 0.01 {\sin\left(\frac{\pi y}{H}\right)}/{1 + \left(\frac{x - x_c}{a}\right)^2}
\end{equation}
with \(x_c = \SI{100}{\kilo\meter}, a = \SI{5}{\kilo\meter}\) and \(H = \SI{10}{\kilo\meter}\). For what concerns the boundary conditions, we consider periodic boundary conditions for the horizontal direction and wall boundary conditions for the vertical direction. A background horizontal velocity \(u = \SI{20}{\meter\per\second}\) is imposed. The grid is composed by \(300 \times 10\) elements with \(r = 4\), while the time step is taken equal to \(\SI{3}{\second}\) yielding \(C \approx 4.17\) and \(C_{u} \approx 0.24\). Figure \ref{fig:IGW}a shows the the potential temperature perturbation at \(t = \SI{3000}{\second}\), where one can easily notice that inertia-gravity waves propagate from the initial perturbation. The results compare well with those available in the literature, see e.g. \cite{melvin:2019, tumolo:2015}. Figure \ref{fig:IGW}b shows the one-dimensional profile of the potential temperature perturbation along \(z = \SI{5}{\kilo\meter}\), which is symmetric about the position \(x = \SI{160}{\kilo\meter}\) and in excellent agreement with the results reported in \cite{giraldo:2008}.  

\subsection{2D density current}
\label{ssec:straka}

\begin{figure}[h!]
	\begin{subfigure}{0.475\textwidth}
		\centering
		\includegraphics[width=0.75\textwidth]{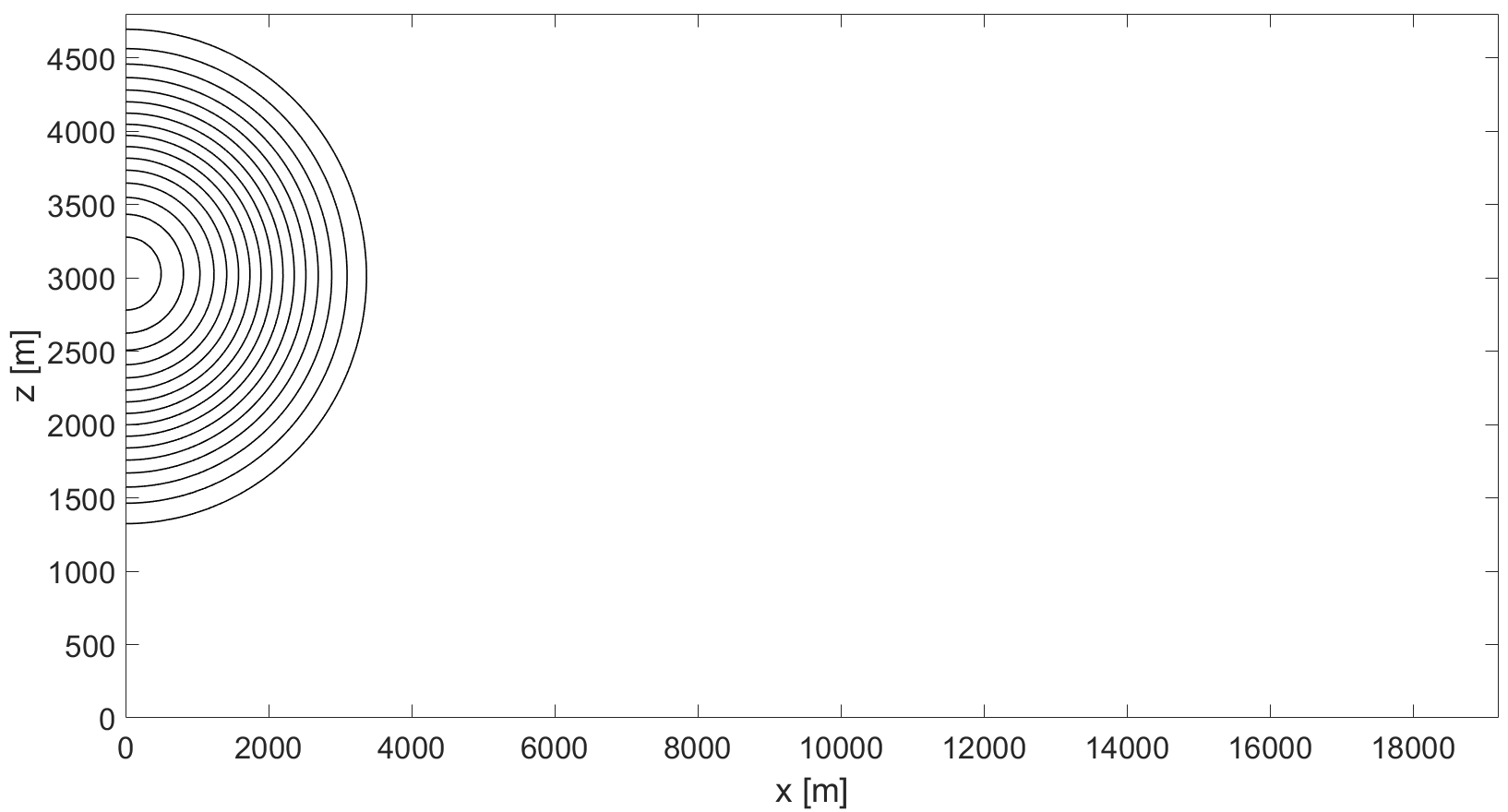} a)
	\end{subfigure}	
	\begin{subfigure}{0.475\textwidth}
		\centering
		\includegraphics[width=0.75\textwidth]{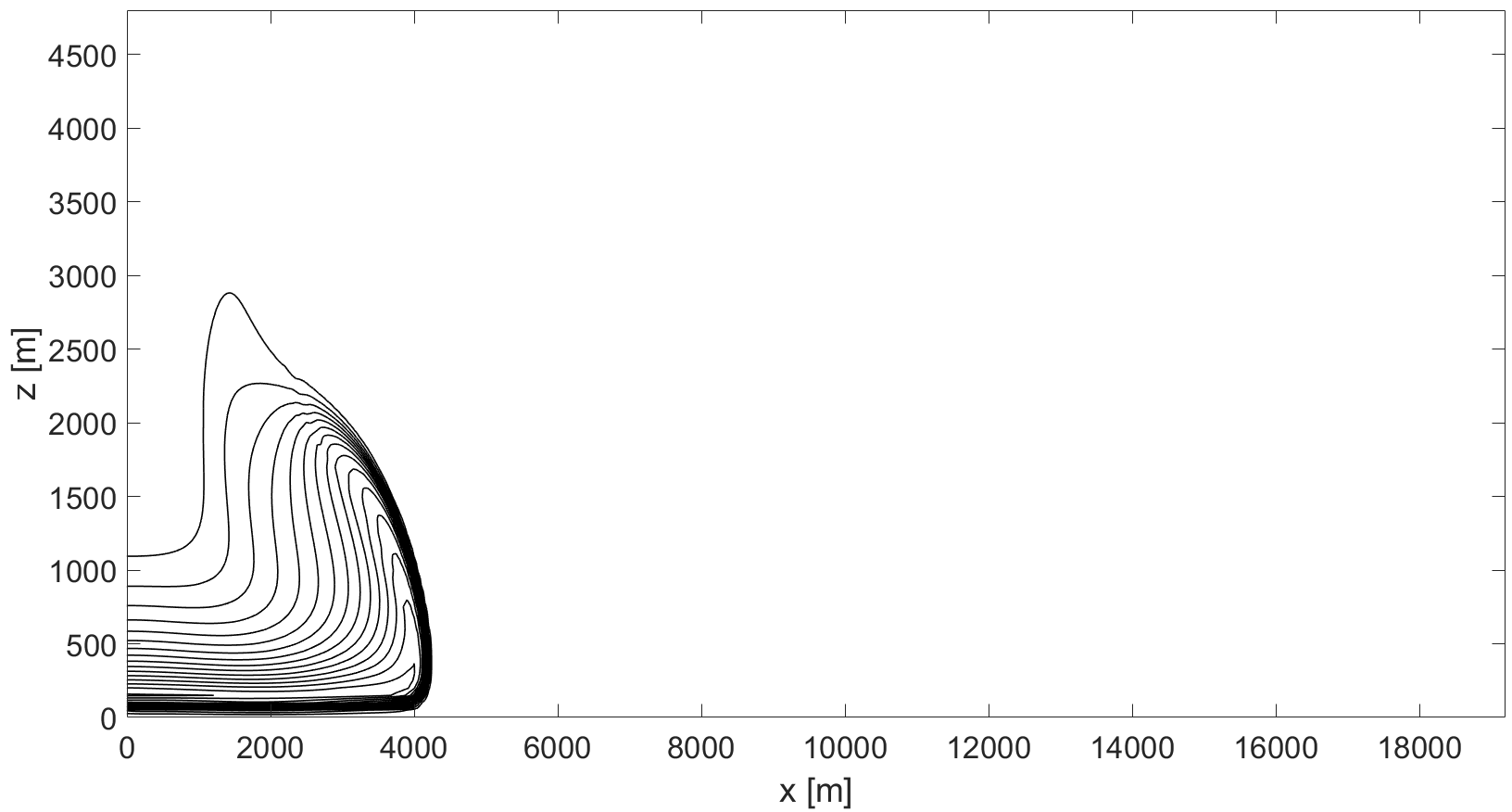} b)
	\end{subfigure}	
	\begin{subfigure}{0.475\textwidth}
		\centering
		\includegraphics[width=0.75\textwidth]{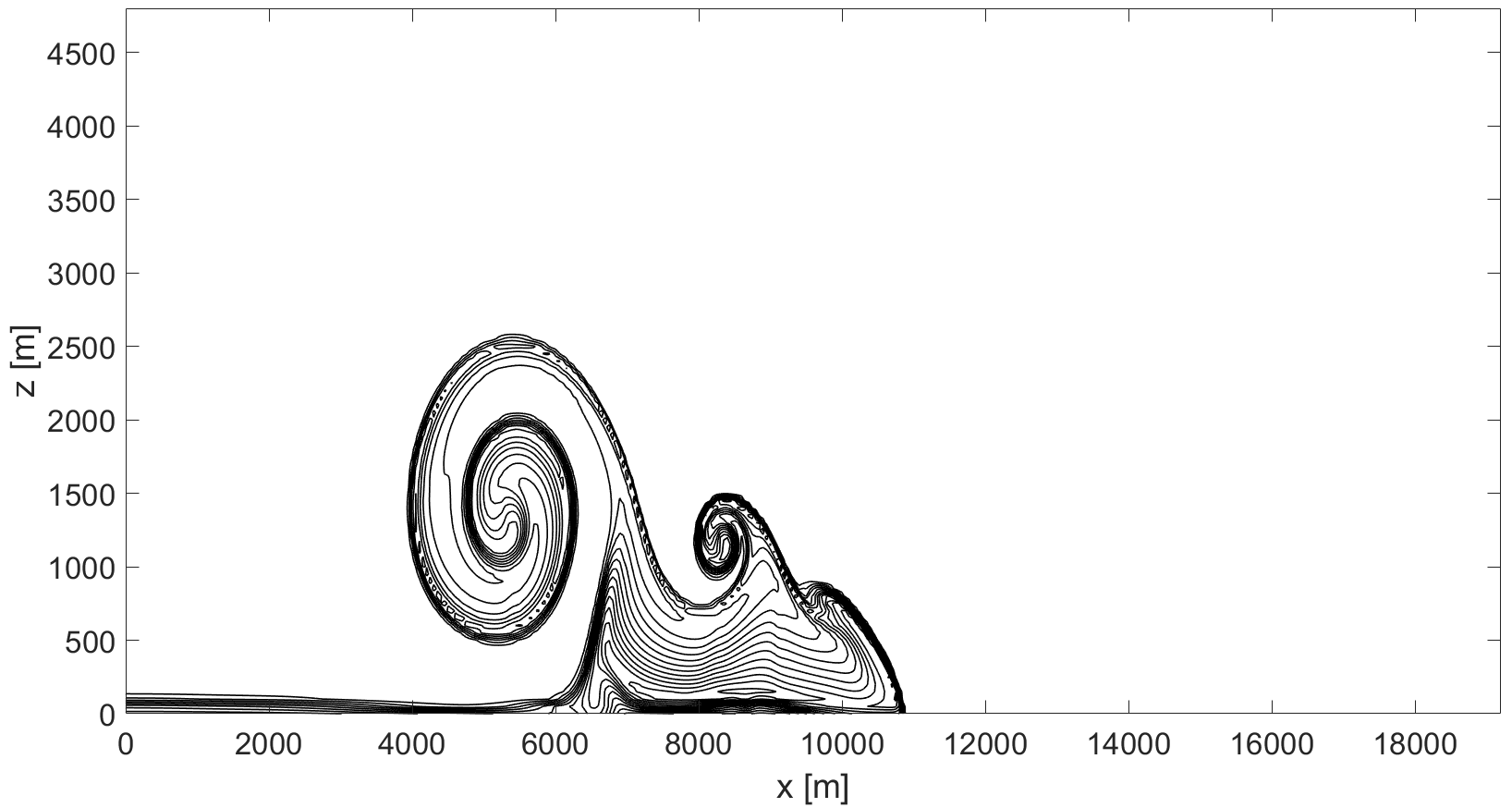} c)
	\end{subfigure}	
	\begin{subfigure}{0.475\textwidth}
		\centering
		\includegraphics[width=0.75\textwidth]{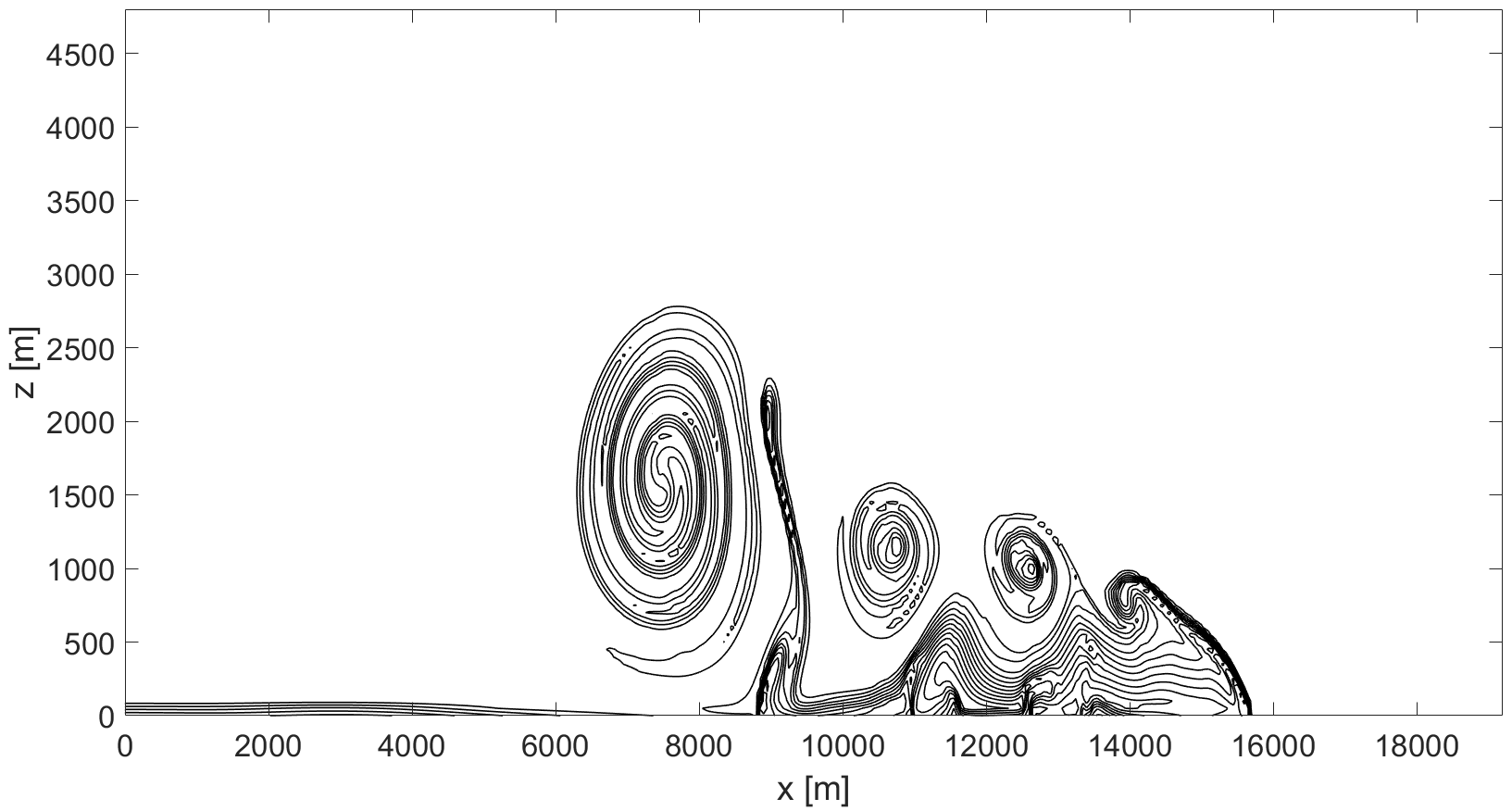} d)
	\end{subfigure}	
	\caption{Current density benchmark, potential temperature deviation from background at: a) \(t = \SI{0}{\second}\), b) \(t = \SI{300}{\second}\), c) \(t = \SI{600}{\second}\), d) \(t = \SI{900}{\second}\). Contours are plotted from \(\SI{-16}{\kelvin}\) to \(\SI{-1}{\kelvin}\) with interval equal to \(\SI{1}{\kelvin}\).}
	\label{fig:Straka_theta_deviation}
\end{figure}

In this Section, we consider the classical density current benchmark proposed in \cite{straka:1993}. The setup consists of a negative temperature perturbation in a motionless isentropic atmosphere with background potential temperature \(\bar{\theta} = \SI{300}{\kelvin}\) and temperature 
$ \bar{T} = \left(300 - zg\frac{\gamma R}{\gamma - 1}\right) \hspace{0.05cm} \SI{}{\kelvin} $
on the domain \(\Omega = \left(-25.6, 25.6\right) \times \left(0, 6.4\right) \SI{}{\kilo\meter}\). The temperature perturbation \(T^{'}\) is defined as
\begin{equation}
T^{'} = 0 \qquad \text{if } \tilde{r} > \SI{1000}{\meter} \ \ \ \ \ \ \
T^{'} =-15\frac{1 + \cos\left(\pi \tilde{r}\right)}{2} \qquad \text{if } \tilde{r} \le \SI{1000}{\meter} 
\end{equation}
where \(\tilde{r} = \sqrt{\left[\frac{\left(x - x_c\right)}{x_r}\right]^2 + \left[\frac{\left(x - y_c\right)}{y_r}\right]^2}\), \(x_c = \SI{0}{\meter}\), \(x_r = \SI{4000}{\meter}\), \(y_c = \SI{3000}{\meter}\) and \(x_r = \SI{2000}{\meter}\). Following \cite{straka:1993}, diffusive terms are included to stabilize the flow, so that the compressible Navier-Stokes equations in conservative form are considered. We resort to an operator splitting approach between the hyperbolic part and the diffusive terms, whose discretization is carried out by the implicit part of the IMEX scheme. This choice is commonly made in numerical methods for atmospherics dynamics, see e.g. \cite{kuehnlein:2019, steppeler:2003}. 	
We consider the diffusion coefficient \(\nu = \SI{75}{\meter^2/\second}\) and we set the thermal conductivity value \(\kappa\) so that the Prandtl number is \(Pr = 0.76\). The boundary conditions are periodic on the left and right boundaries and wall boundary conditions on the top and bottom boundaries. The grid is composed by \(1024 \times 128\) elements with \(r = 2\) leading to a resolution equal to \(\SI{25}{\meter}\). The time step is taken equal to \(\SI{0.1}{\second}\), yielding a maximum Courant number \(C \approx 1.4\) and \(C_u \approx 0.15\). 
Figure \ref{fig:Straka_theta_deviation} shows the deviation of the potential temperature with respect to the background value at different times for the subdomain \(\left(0, 19.2\right) \times \left(0, 4.8\right)\SI{}{\kilo\meter}\). In view of the negative buoyancy, the structure falls and reaches the bottom boundary. It then moves to the right, developing vortices. The front location is located at \(x = \SI{15700}{\meter}\), in agreement with the results obtained in \cite{benacchio:2014, melvin:2019}.

\begin{figure}[h!]
	\begin{subfigure}{0.475\textwidth}
		\centering
		\includegraphics[width=0.9\textwidth]{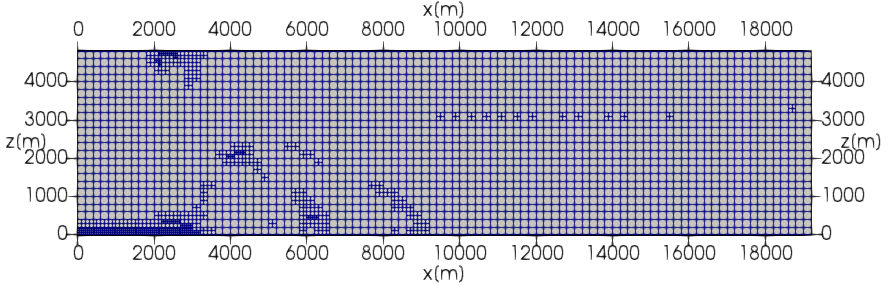} a)
	\end{subfigure}	
	\begin{subfigure}{0.475\textwidth}
		\centering
		\includegraphics[width=0.9\textwidth]{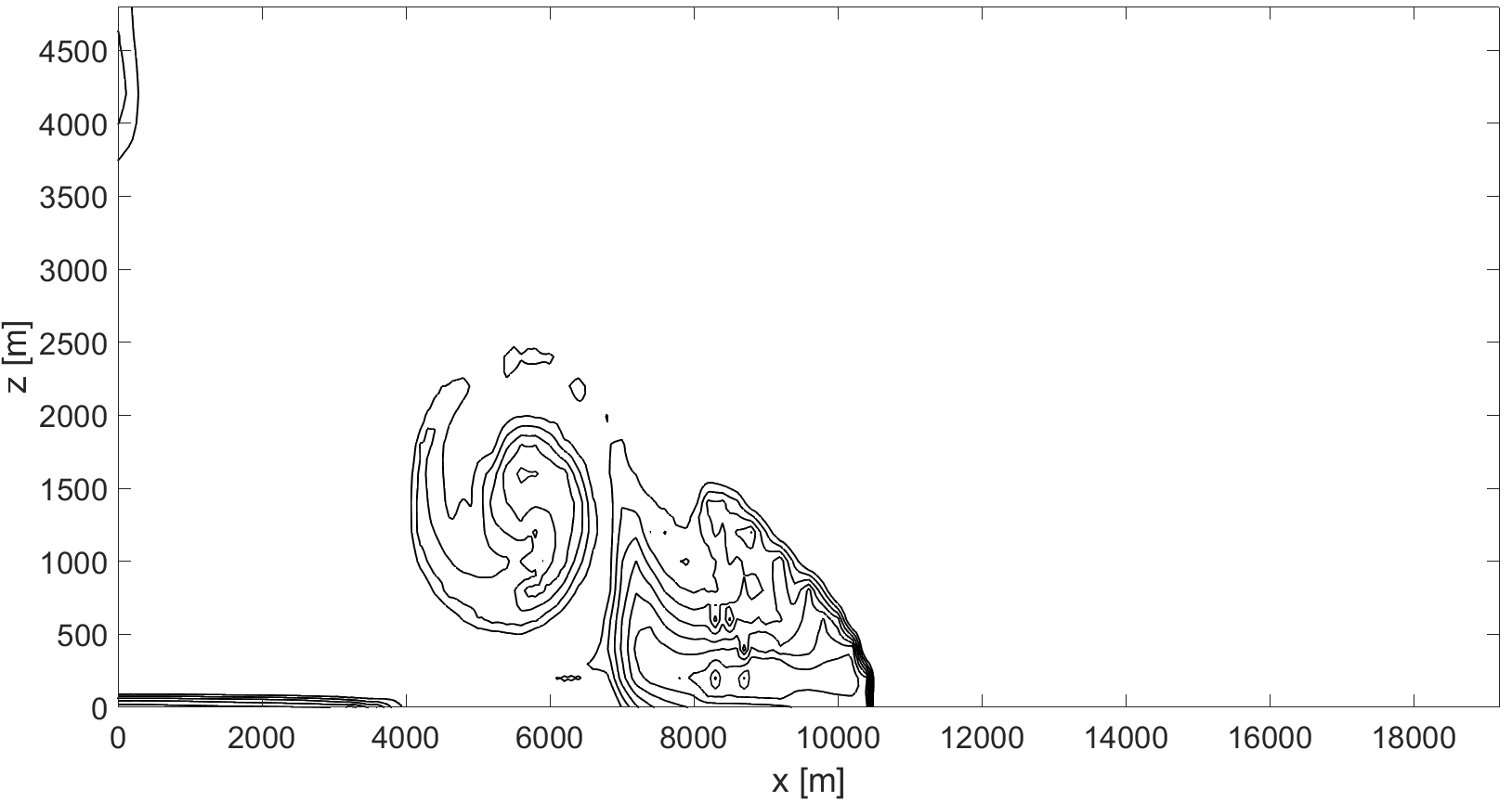} b)
	\end{subfigure}	
	\begin{subfigure}{0.475\textwidth}
		\centering
		\includegraphics[width=0.9\textwidth]{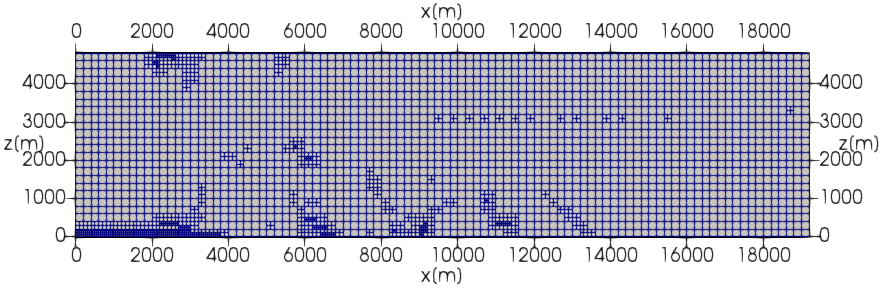} c)
	\end{subfigure}	
	\begin{subfigure}{0.475\textwidth}
		\centering
		\includegraphics[width=0.9\textwidth]{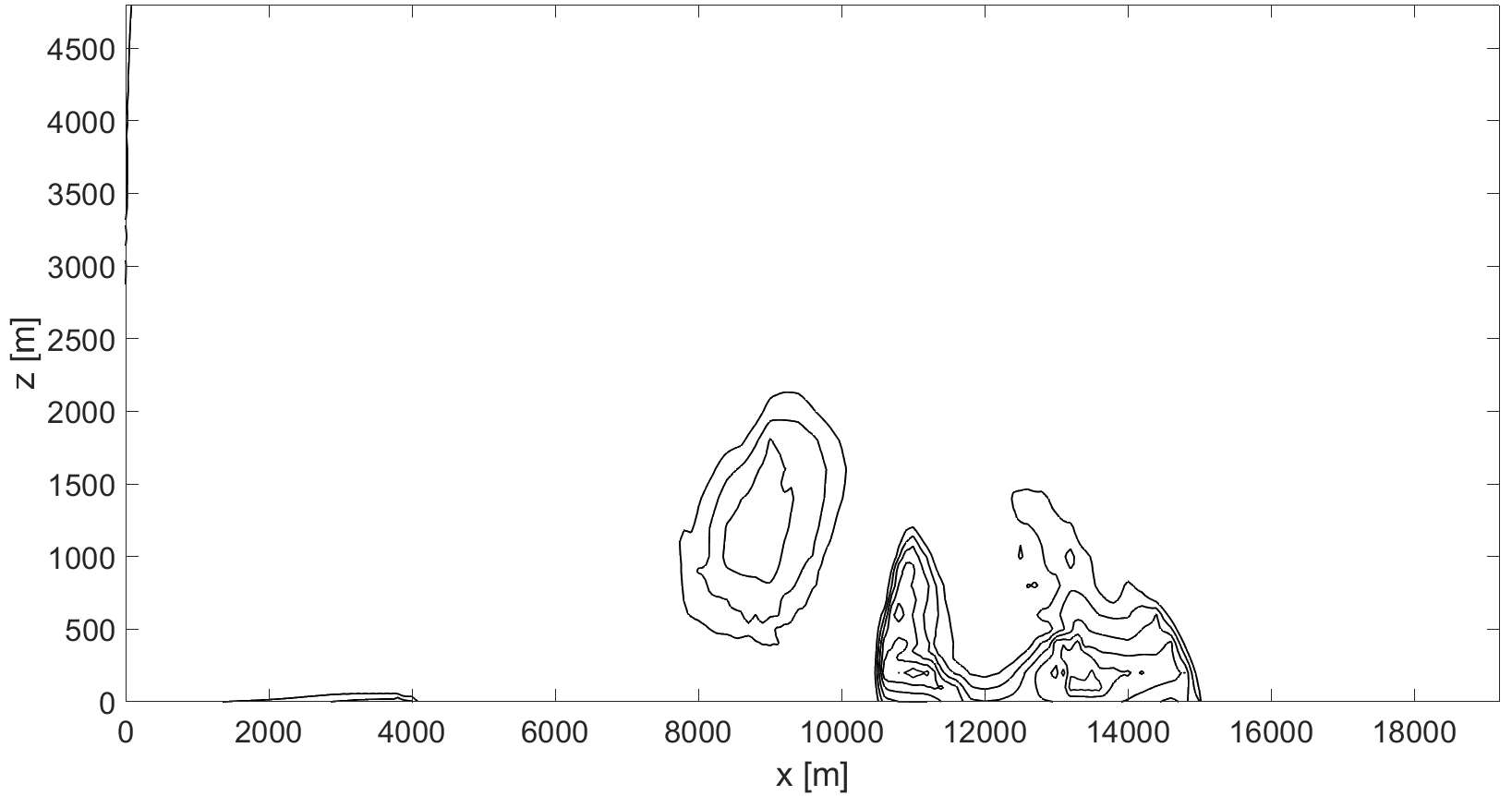} d)
	\end{subfigure}	
	\caption{Current density benchmark; on the left the adaptive meshes at \(t = \SI{600}{\second}\) (a) and \(t = \SI{900}{\second}\) (c) are reported, whereas on the right potential temperature deviations from the background at \(t = \SI{600}{\second}\) (b) and \(t = \SI{900}{\second}\) (d) are reported.}
	\label{fig:Straka_adaptive}
\end{figure}

\begin{figure}[h!]
	\centering
	\includegraphics[width=0.9\textwidth]{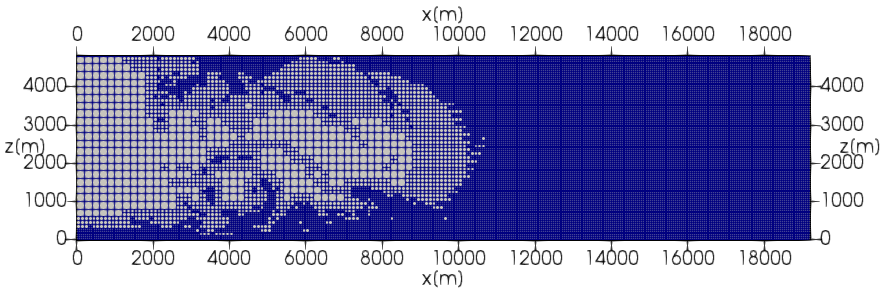} a)
	\caption{Current density benchmark, adaptive grid at \(t = \SI{600}{\second}\) obtained when $\eta_{K}$ exceeds $4 \cdot 10^{-2} \SI{}{\joule\per\meter}$.}
	\label{fig:Straka_adaptive_bis}
\end{figure}

We now demonstrate the $h-$adaptivity capabilities available in the proposed implementation. We use as refinement indicator the gradient of the potential temperature. More specifically, we set
\begin{equation}\label{eq:eta_K}
\eta_K = \max_{i \in \mathcal{N}_{K}} \left|\nabla\theta\right|_{i} 
\end{equation}
where \(\mathcal{N}_{K}\) denotes the set of nodes over the element \(K\). We allow to refine when $\eta_K$ exceeds $3 \cdot 10^{-2} \SI{}{\joule\per\meter}$ and to coarsen when the indicator is below $\cdot 10^{-2} \SI{}{\joule\per\meter}$. The initial computational grid is composed by \(256 \times 64\) elements and only two local refinements are allowed, so as to keep under control the advective Courant number. Figure \ref{fig:Straka_adaptive} shows the adaptive mesh and the contour plot of the potential temperature perturbation at \(t = \SI{600}{\second}\) and at \(t = \SI{900}{\second}\). Figure \ref{fig:Straka_adaptive} shows the adaptive grid and the contour plot of the potential temperature perturbation at \(t = \SI{600}{\second}\) and at \(t = \SI{900}{\second}\). One can notice that more resolution is added in correspondence of the bubble. A computational saving time of around \(60\%\) is achieved by the adaptive mesh refinement procedure. It is worth to notice that, for this test case, the mesh adaptation procedure is strongly dependent on the bounds for the refinement criterion. Indeed, allowing to refine when $\eta_K$ exceeds $4 \cdot 10^{-2} \SI{}{\joule\per\meter}$ leads to an excess of resolution in the domain, as evident from Figure \ref{fig:Straka_adaptive_bis}. 

\subsection{2D falling bubble}
\label{ssec:cold_bubble}

\begin{figure}[h!]
	\begin{subfigure}{0.475\textwidth}
		\centering
		\includegraphics[width=0.75\textwidth]{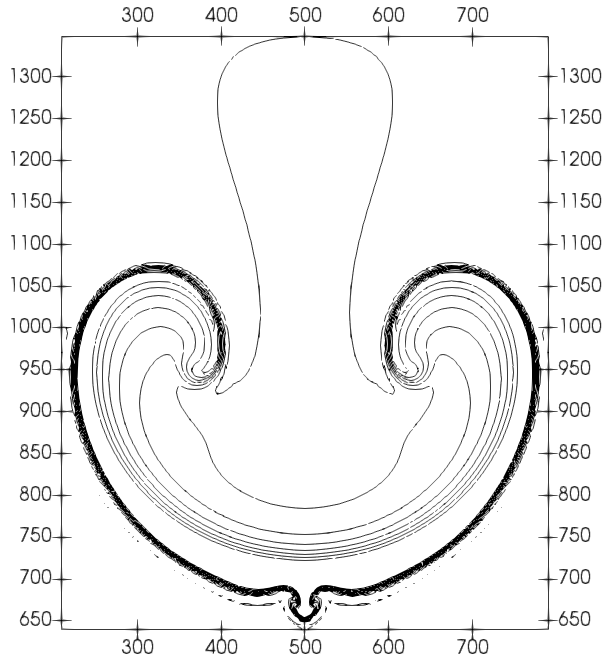} a)
	\end{subfigure}	
	\begin{subfigure}{0.475\textwidth}
		\centering
		\includegraphics[width=0.75\textwidth]{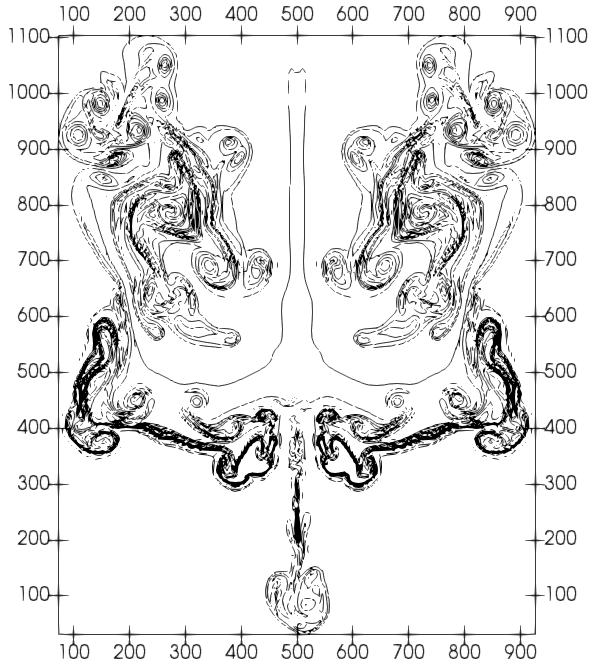} b)
	\end{subfigure}	
	\caption{Cold bubble benchmark, potential temperature deviation from background at:  a) \(t = \SI{100}{\second}\), b) \(t = \SI{200}{\second}\). Contours are plotted  from \(\SI{-11.5}{\kelvin}\) to \(\SI{8}{\kelvin}\) with interval equal to \(\SI{1.625}{\kelvin}\).}
	\label{fig:Cold_Bubble_theta_deviation}
\end{figure}

In this Section, we consider a test case proposed in \cite{restelli:2009}. The computational domain is the rectangle \(\left(0,1000\right) \times \left(0,2000\right) \SI{}{\meter}\) and the initial condition is represented by a thermal anomaly introduced in an isentropic background atmosphere with constant potential temperature \(\bar{\theta} = \SI{303}{\kelvin}\). The perturbation of potential temperature \(\theta^{'}\) defines the initial datum and it is given by
\begin{equation}
\theta^{'} = 
A  \ \ \text{if } \tilde{r} \le r_0 \ \ \ \ 
\theta^{'} = A\exp\left(-\frac{\left(\tilde{r} - r_0\right)^2}{\sigma^2}\right) \ \  \text{if } \ \ \tilde{r} > r_0,
\end{equation}
with \(\tilde{r}^2 = \left(x - x_0\right)^2 + \left(z - z_0\right)^2\) and \(x_0 = \SI{500}{\meter}\), \(z_0 = \SI{1250}{\meter}\), \(r_0 = \SI{50}{\meter}\), \(\sigma = \SI{100}{\meter}\) and \(A = \SI{-15}{\kelvin}\). The expression of the initial profile of the Exner pressure is given by 
$ \overline{\pi} = 1 - \frac{g}{c_{p}\theta}z, $
with \(c_{p} = \frac{\gamma}{\gamma - 1}R = \SI{1004.5}{\joule\per\kilogram\per\kelvin}\) denoting the specific heat at constant pressure. Notice that, unlike in \cite{restelli:2009}, no artificial viscosity has been added to stabilize the computation. Wall boundary conditions are imposed at all the boundaries. The time step is taken to be \(\Delta t = \SI{0.08}{\second}\), corresponding to a maximum Courant number \(C \approx 5.6\) and a maximum advective Courant number \(C_u \approx 0.24\), whereas the final time is \(T_{f} = \SI{200}{\second}\). The computational grid is composed by \(200 \times 400\) elements with \(r = 1\) leading to a resolution equal to \(\SI{5}{\meter}\). Figure \ref{fig:Cold_Bubble_theta_deviation} shows the contours of the potential temperature deviation from the background at \(t = \SI{0}{\second}\), \(t = \SI{100}{\second}\) and \(t = \SI{200}{\second}\) and the results are in reasonable agreement with those reported in \cite{restelli:2009}. 

\begin{figure}[h!]
	\begin{subfigure}{0.475\textwidth}
		\centering
		\includegraphics[width=0.75\textwidth]{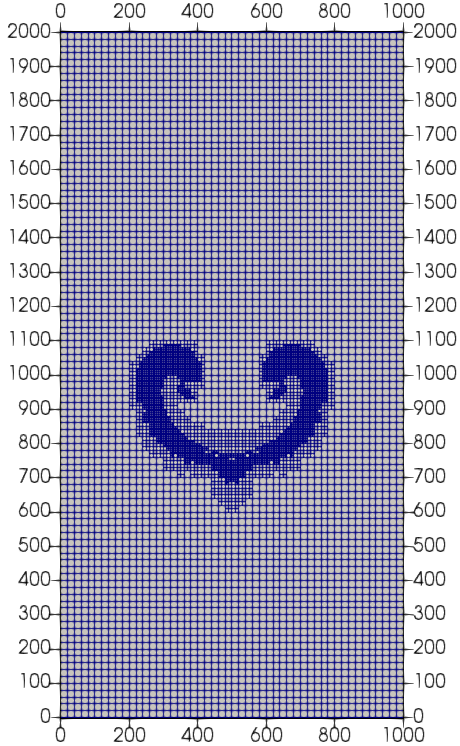} a)
	\end{subfigure}	
	\begin{subfigure}{0.475\textwidth}
		\centering
		\includegraphics[width=0.75\textwidth]{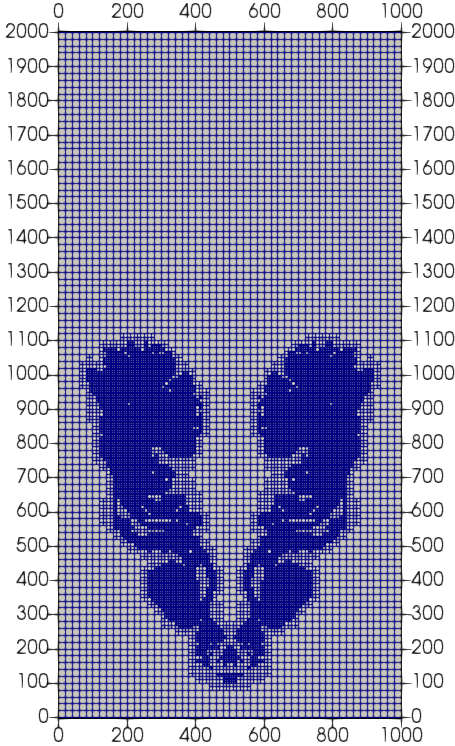} b)
	\end{subfigure}	
	\caption{Cold bubble benchmark, adaptive mesh at: a) \(t = \SI{100}{\second}\), b) \(t = \SI{200}{\second}\).}
	\label{fig:Cold_Bubble_mesh_adaptive}
\end{figure}

\begin{figure}[h!]
	\includegraphics[width=0.30\textwidth]{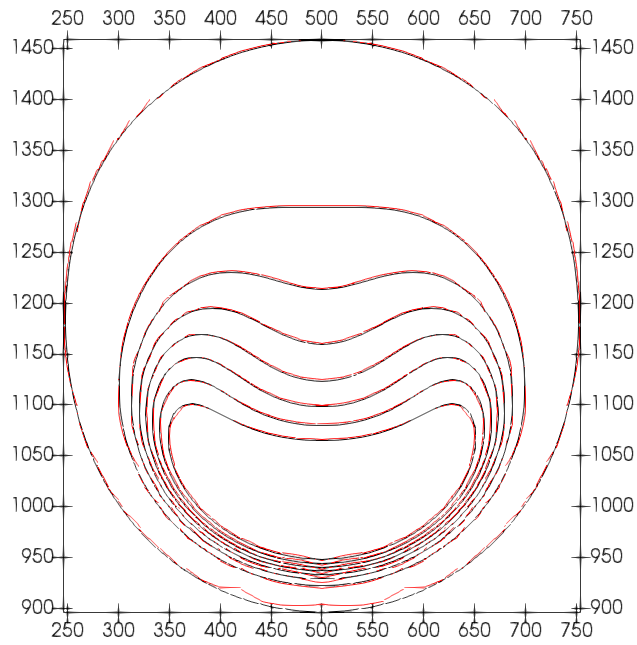} a)
	\includegraphics[width=0.28\textwidth]{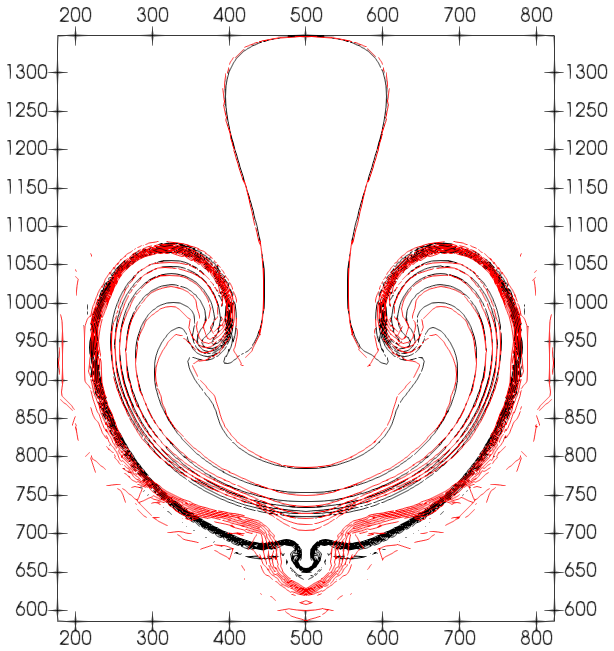} b)
	\includegraphics[width=0.26\textwidth]{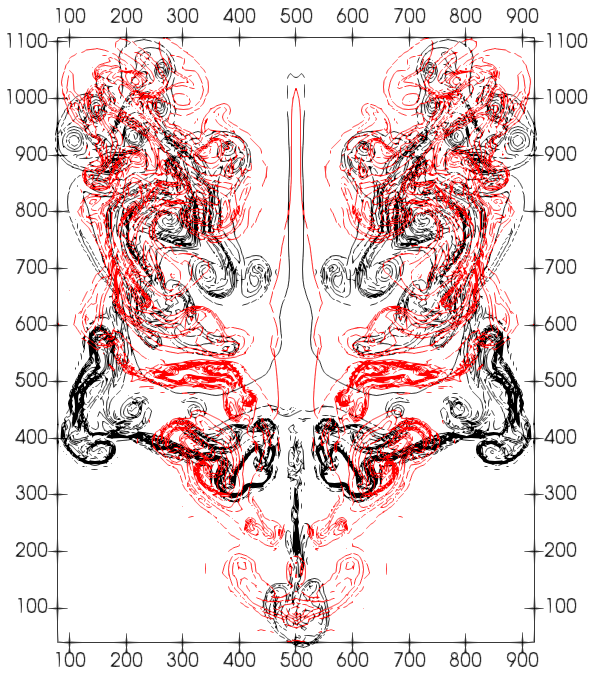} c)
	\caption{Cold bubble benchmark, potential temperature deviation from background at: a) \(t = \SI{50}{\second}\), b) \(t = \SI{100}{\second}\), c) \(t = \SI{200}{\second}\). Black line denotes the results with uniform grid, whereas red line represents the results with \(h\)-adaptivity. Contours are plotted from \(-11.5\) to \(8\) with interval equal to \(\SI{1.625}{\kelvin}\).}
	\label{fig:Cold_Bubble_theta_deviation_adaptive}
\end{figure}

We consider now the $h-$adaptivity version of the proposed scheme. We use the same refinement indicator \eqref{eq:eta_K} introduced in Section \ref{ssec:straka}. We allow to refine when $\eta_K$ exceeds $10^{-4} \hspace{0.05cm} \SI{}{\joule\per\meter}$ and to coarsen when the indicator is below $6 \cdot 10^{-5} \hspace{0.05cm} \SI{}{\joule\per\meter}$. The initial computational grid is composed by \(50 \times 100\) elements and only two local refinements are allowed, so as to control the advective Courant number and to match at the finest refinement level the resolution of the non adaptive mesh simulation. As one can easily notice from Figure \ref{fig:Cold_Bubble_mesh_adaptive}, the refinement criterion is able to track the bubble. The contour plots in Figure \ref{fig:Cold_Bubble_theta_deviation_adaptive} show a reasonable agreement for \(t = \SI{50}{\second}\) and \(t = \SI{100}{\second}\), whereas for \(t = \SI{200}{\second}\) significant differences between the simulations with uniform and adaptive grid appear, which are due to the different development of the Kelvin-Helmholtz instability. The final mesh consists of 19334 elements instead of the 80000 elements of the full resolution mesh and achieves a computational time reduction of 25\%.

\subsection{3D rising bubble}
\label{ssec:3D_bubble}

\begin{figure}[h!]
	\begin{subfigure}{0.475\textwidth}
		\centering
		\includegraphics[width=0.72\textwidth]{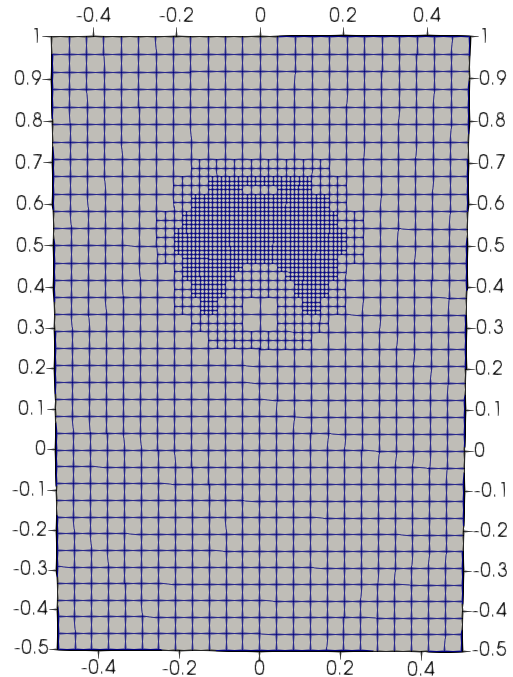} a)
	\end{subfigure}	
	\begin{subfigure}{0.475\textwidth}
		\centering
		\includegraphics[width=0.82\textwidth]{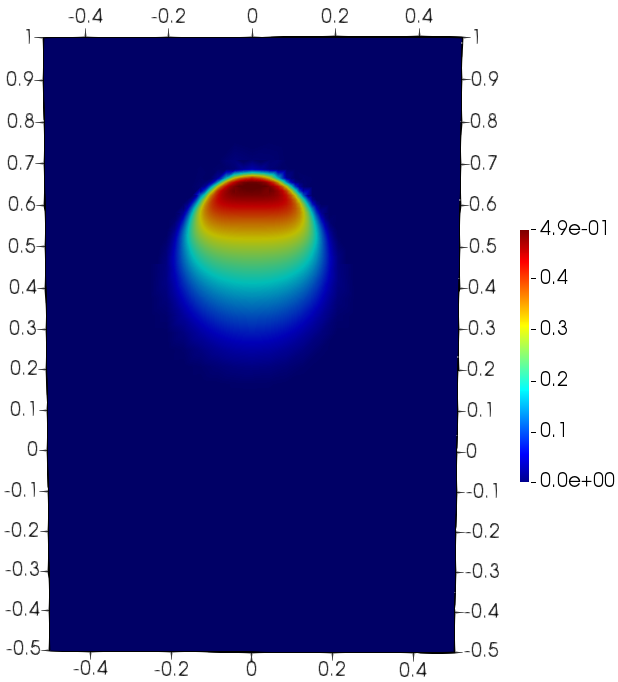} b)
	\end{subfigure}	
	\begin{subfigure}{0.475\textwidth}
		\centering
		\includegraphics[width=0.7\textwidth]{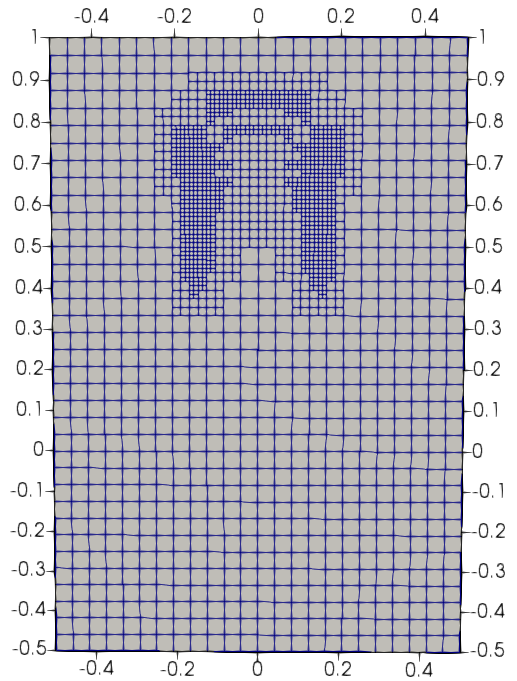} c)
	\end{subfigure}	
	\begin{subfigure}{0.475\textwidth}
		\centering
		\includegraphics[width=0.8\textwidth]{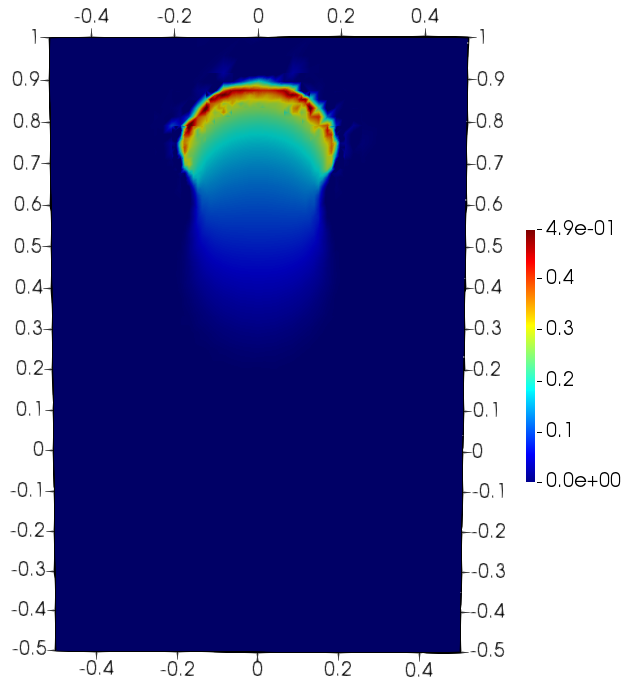} d)
	\end{subfigure}	
	\caption{3D rising bubble benchmark, results for \(y = \SI{0}{\meter}\); on the left the adaptive meshes at \(t = \SI{200}{\second}\) (a) and \(t = \SI{400}{\second}\) (c) are reported, whereas on the right potential temperature deviations from the background at \(t = \SI{200}{\second}\) (b) and \(t = \SI{400}{\second}\) (d) are reported.}
	\label{fig:3D_Bubble_contours}
\end{figure}

In this Section, we consider the 3D rising bubble benchmark proposed in \cite{melvin:2019}. A neutrally stratified isentropic atmosphere is assumed, with \(\bar{\theta} = \SI{300}{\kelvin}\) in the domain \(\Omega = \left(-500, -500, -500\right) \times \left(500, 500, 1000\right) \SI{}{\meter}\). A spherical perturbation \(\theta^{'}\) located at \(\left(x_0, y_0, z_0\right) = \left(0, 0, 350\right)\SI{}{\meter}\) is added to the potential temperature
\begin{equation}
\theta^{'} = 
0.25\left[ 1 + \cos\left(\frac{\pi r}{r_0}\right) \right] \qquad \text{if } \ \ \tilde{r} \le r_0 \ \ \ \ 
\theta^{'} = 0 \qquad \text{if }  \ \ \tilde{r} > r_0,
\end{equation}
with \(\tilde{r} = \sqrt{\left(x - x_0\right)^2 + \left(y - y_0\right)^2 + \left(z - z_0\right)^2}\) and \(r_0 = \SI{250}{\meter}\). Wall boundary conditions are imposed for all the six boundaries and we take \(r = 2\). In order to enhance the computational efficiency, we use the \(h-\)adaptivity capabilities, with the same refinement indicator \eqref{eq:eta_K} introduced in Section \ref{ssec:straka}. 
For this benchmark, we allow to refine when $\eta_K$ exceeds $3 \cdot 10^{-3}$ and to coarsen when the indicator is below $5 \cdot 10^{-4}$. The initial grid is composed by \(24 \times 24 \times 36\) elements and we allowed up to two local refinements, whose finest refinement level corresponds to a uniform mesh with \(96 \times 96 \times 144\) elements and to a resolution around \(\SI{5}{\meter}\). The time step is taken to be equal to \(\Delta t=\SI{0.4}{\second}\), leading to a maximum acoustic Courant number \(C \approx 27\) and advective Courant number \(C_u \approx 0.22\). Figure \ref{fig:3D_Bubble_contours} shows snapshots of the bubble at \(t = \SI{200}{\second}\) and \(t = \SI{400}{\second}\). At the later time, a Kelvin-Helmholtz instability starts to develop, which is however still insufficiently well resolved by the present mesh. Further refinement levels 
or higher polynomial degrees will have to be employed in future simulations to achieve better accuracy at the later stage. Similar issue for an analogous test case are reported in \cite{busto:2020}. The final grid is composed by \(62792\) elements. 

The size of this benchmark makes it a good candidate for a parallel scaling test. An initial mesh composed by \(48 \times 48 \times 72\) elements corresponding to \(13436928\) dofs for the velocity and \(4478976\) dofs for the remaining scalar variables is considered. Two configurations are employed: in the first case we keep it fixed, whereas in the second one we apply \(h-\)adaptivity with two local refinements, roughly doubling the number of degrees of freedom. A strong scaling analysis is performed executing the simulation up to time \(t = \SI{8}{\second}\) and we use from 32 up to 1024 2xCPU x86 Intel Xeon Platinum 8276-8276L @ 2.4Ghz cores of the HPC infrastructure GALILEO100 at the Italian supercomputing center CINECA. 

\begin{figure}[h!]
	\centering
	\includegraphics[width=0.7\textwidth]{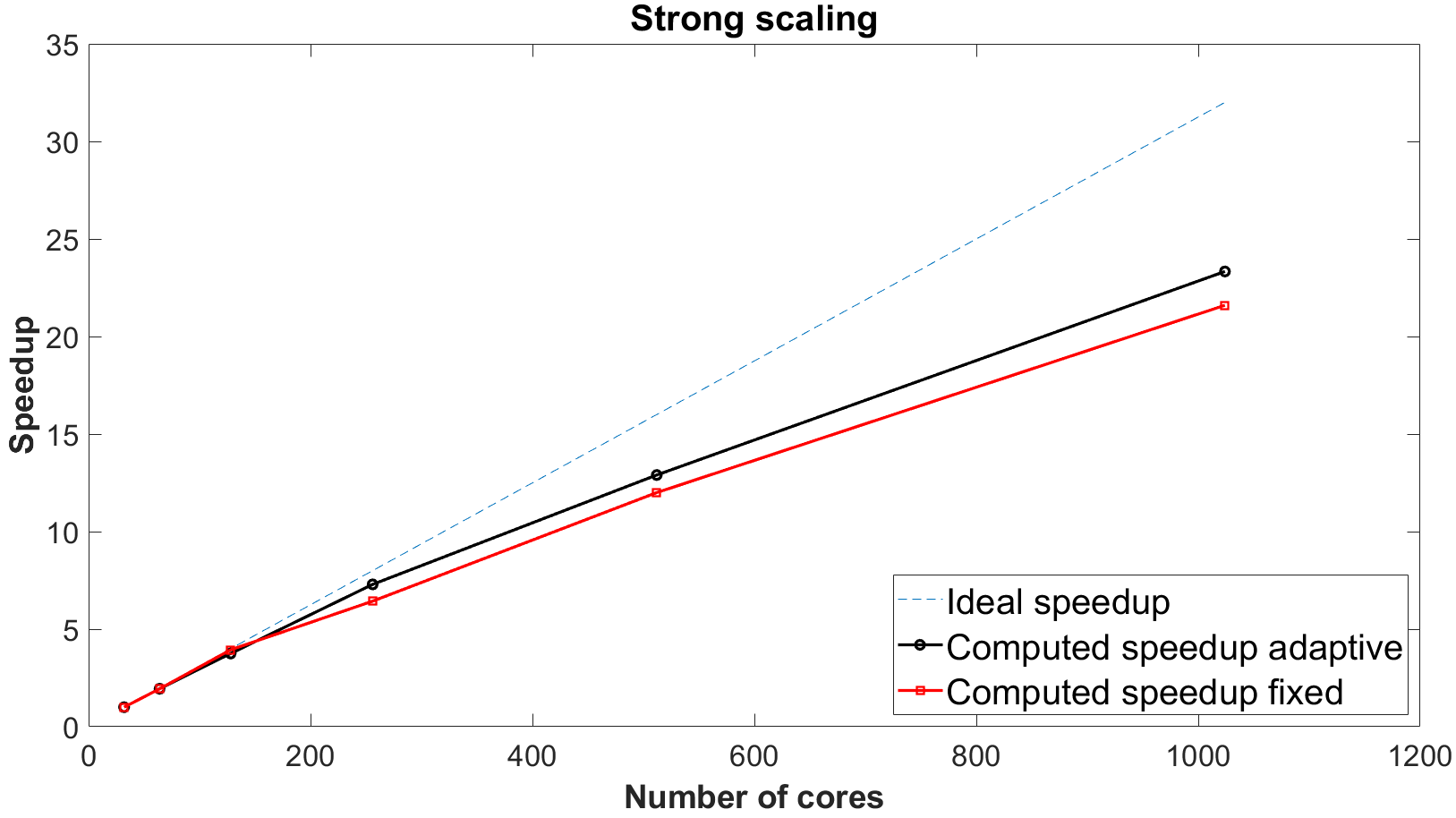} 
	\caption{3D rising bubble benchmark, strong scaling analysis. The speedup is computed with respect to the time required with \(32\) cores.}
	\label{fig:strong_scaling}
\end{figure}

The results, reported in Figure \ref{fig:strong_scaling} and Table \ref{tab:strong_scaling_times}, are quite similar for the two configurations. A good linear scaling is obtained up to 128 cores, even with superlinear behaviour for the fixed mesh framework due to cache effects. Starting from 256 cores, the performance of the fixed mesh configuration exhibits a small degradation and, for a higher number of cores, the speed-up is less optimal for both configurations, due to overwhelming communication costs. The apparent better behaviour of the \(h-\)adaptive version is due to the fact that more degrees of freedoms are involved and, therefore, the role of communication costs is less evident. The result also highlights that the local refining procedure has no significant impact on the parallel performance and that both efficiency and scalability can be achieved in this framework. 

\begin{table}[h!]
	\centering
	\begin{tabular}{|c|c|c|}
		\hline
		Number of cores & Wallclock time fixed grid  & Wallclock time adaptive grid \\
		\hline
		32 & \(5.40 \times 10^{3}\) & \(3.29 \times 10^{3}\) \\
		\hline
		64 & \(2.78 \times 10^{3}\) & \(1.69 \times 10^{3}\) \\
		\hline
		128 & \(1.37 \times 10^{3}\) &  \(8.77 \times 10^{2}\) \\
		\hline
		256 & \(8.39 \times 10^{2}\) & \(4.51 \times 10^{2}\) \\
		\hline
		512 & \(4.50 \times 10^{2}\) & \(2.55 \times 10^{2}\) \\
		\hline
		1024 & \(2.50 \times 10^{2}\) & \(1.41 \times 10^{2}\) \\
		\hline
	\end{tabular}
	\caption{Wallclock times in seconds for the strong scaling analysis.}
	\label{tab:strong_scaling_times}
\end{table}

\section{Numerical results with orography}
\label{sec:tests_with} \indent

We now consider a number of tests concerning idealized flows over orography, that since the seminal papers \cite{klemp:1983, klemp:1978} have become a standard benchmark for numerical models of atmospheric flows, see e.g. the results and discussions in \cite{bonaventura:2000, melvin:2019, pinty:1995, tumolo:2015}.
In most of the tests, the bottom boundary is described by the function
\begin{equation}\label{eq:versiera_Agnesi}
h(x) =  \frac{h_{m}}{1 + \left(\frac{x - x_{c}}{a_{c}}\right)^2},
\end{equation}
the so-called \textit{versiera di Agnesi}, where \(h_{m}\) is the height of the hill and \(a_{c}\) is the half-width. The classical Gal-Chen height-based terrain-following coordinate \cite{galchen:1975} is used to build the mapping between the reference element and the physical one and to obtain a terrain following mesh in Cartesian coordinates.

\subsection{2D hydrostatic flow over a hill}
\label{ssec:hydrostatic}

\begin{figure}[h!]
	\centering
	\includegraphics[width=0.7\textwidth]{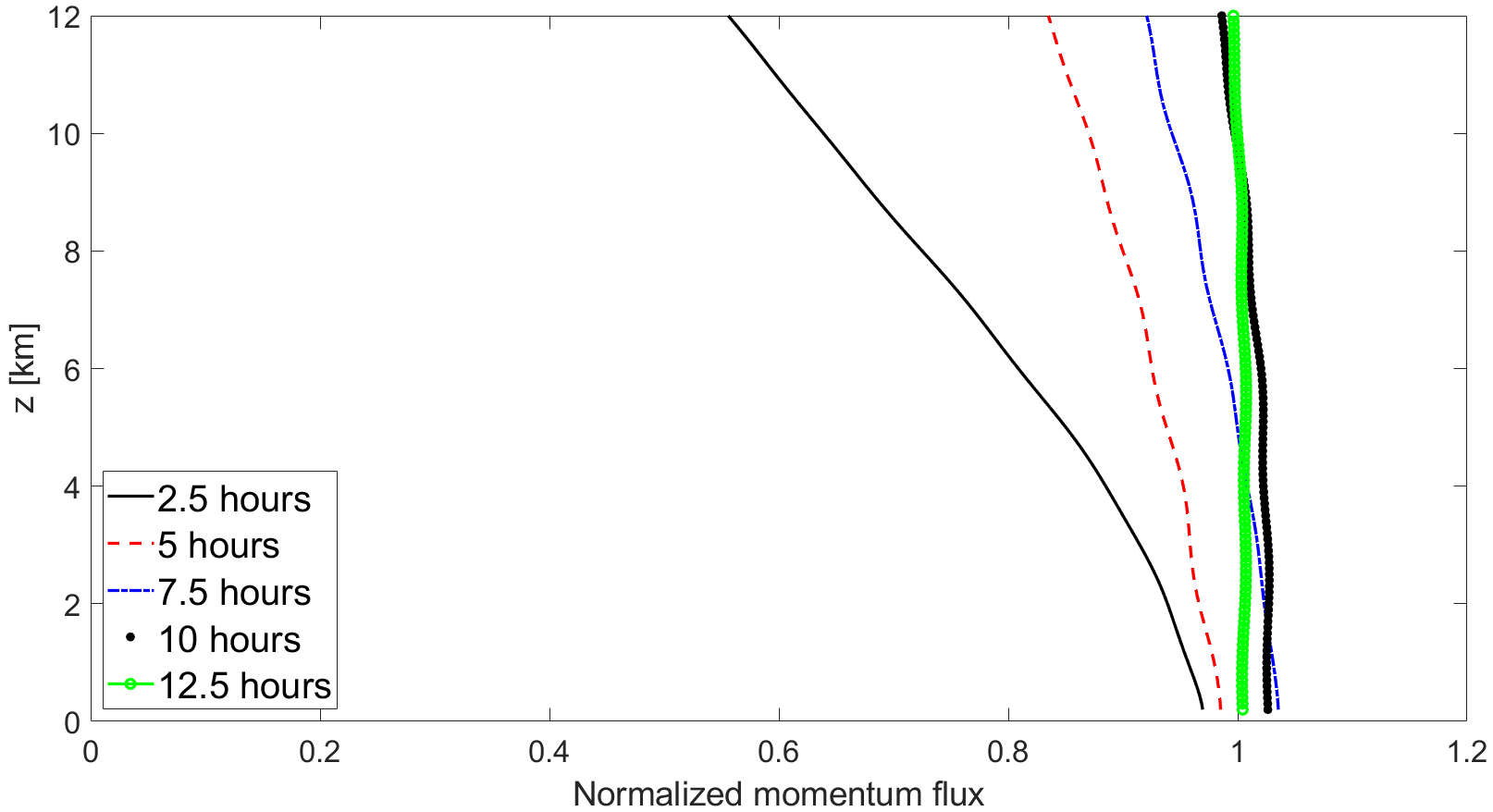}
	\caption{Linear hydrostatic flow over a hill, evolution of normalized momentum flux.}
	\label{fig:linear_hydro_momentum}
\end{figure}

We first consider the linear hydrostatic configuration presented e.g. in \cite{giraldo:2008}. The computational domain is \(\Omega = \left(0, 240\right) \times \left(0, 30\right) \hspace{0.05cm} \SI{}{\kilo\meter}\) with \(h_{m} = \SI{1}{\meter}, x_{c} = \SI{120}{\kilo\meter}\) and \(a_{c} = \SI{10}{\kilo\meter}\). The final time is \(T_{f} = \SI{45000}{\second}\). The initial state of the atmosphere consists of a constant mean flow with \(\overline{u} = \SI{20}{\meter\per\second}\) and of an isothermal background profile with temperature \(\overline{T} = \SI{250}{\kelvin}\). The initial profile of the Exner pressure is given by
$ \overline{\pi} = \left( {p_{0}}/{\overline{p}}\right)^{\frac{\gamma - 1}{\gamma}} = \exp\left(-\frac{g}{c_{p}\overline{T}}z\right). $
We recall that \(c_{p} = \frac{\gamma}{\gamma - 1}R\) denotes the specific heat at constant pressure and that here \(p_{0} = 10^{5} \hspace{0.05cm} \SI{}{\pascal}\); moreover, since in an isothermal configuration the Brunt-V{\"{a}}is{\"{a}}l{\"{a}} frequency is \(N = g/\sqrt{c_{p}\overline{T}}\), it can be easily checked that \(\frac{N a_{c}}{\overline{u}} >> 1\), so that this configuration corresponds to a hydrostatic regime according to the classification in \cite{pinty:1995}. For what concerns the boundary conditions, wall boundary conditions are used for the bottom boundary and non-reflecting boundary conditions are required by the top boundary and the lateral boundaries. For this purpose, we introduce a Rayleigh damping profile following \cite{melvin:2019}, so that
\begin{equation}
\lambda = 
0, \qquad  \text{if } \ \ z < z_{B}  \ \ \ \ \ 
\lambda =  \overline{\lambda}\sin^2\left[\frac{\pi}{2}\left(\frac{z - z_{B}}{z - z_{T}}\right)\right] \qquad  \text{if } \ \ z \ge z_{B},
\end{equation}
where \(z_{B}\) denotes the height at which the damping starts and \(z_{T}\) is the top height of the considered domain. Analogous definitions apply for the two lateral boundaries. In this case, we consider \(\overline{\lambda}\Delta t = 0.3\) and we apply the damping layer in the topmost \(\SI{15}{\kilo\meter}\) of the domain and in the first and last \(\SI{80}{\kilo\meter}\) along the horizontal direction. The grid is composed by \(100 \times 75\) elements with \(r = 4\), yielding a resolution of \(\SI{600}{\meter}\) along \(x\) and \(\SI{100}{\meter}\) along $z, $ whereas the time-step is equal to \(\SI{2.5}{\second}\), leading to \(C \approx 1.84\) and \(C_{u} \approx 0.12\). Following \cite{smith:1979}, we also define the vertical momentum flux as
\begin{equation}
m(z) = \int_{-\infty}^{\infty} \overline{\rho}(z)u^{'}(x,z)w^{'}(x,z)dx,
\end{equation} 
where \(u^{'}\) and \(w^{'}\) represent the deviation from the background state of the horizontal and vertical velocity, respectively.
This is a very important diagnostic quantity in atmospheric modelling, used to check that the numerical model is correctly reproducing the orographic forcing on the main flow. From the linear theory, the analytical momentum flux is given by
\begin{equation}\label{eq:analytic_momentum_hydro}
m^{H} = -\frac{\pi}{4}\overline{\rho}_{s}\overline{u}_{s} N h_{m}^{2},
\end{equation}
where \(\overline{\rho}_{s}\) and \(\overline{u}_{s}\) denote the surface background density and velocity, respectively. Figure \ref{fig:linear_hydro_momentum} shows the behaviour over time of the momentum flux normalized by its analytical value. It can be noticed that the analytical value is approached as the simulation reaches the steady state. 

\begin{figure}[h!]
	\centering
	\includegraphics[width=0.45\textwidth]{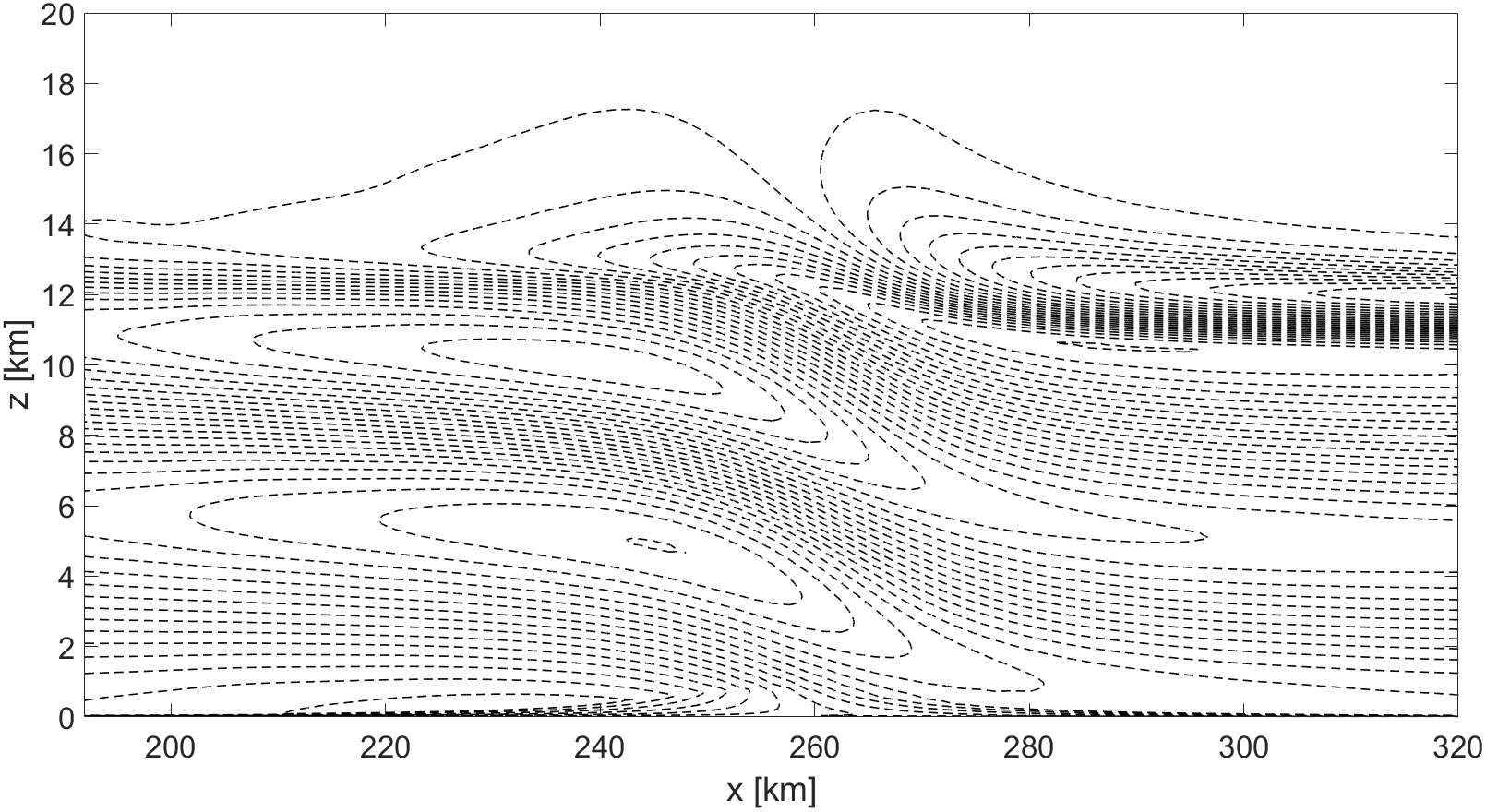} a)
	\includegraphics[width=0.45\textwidth]{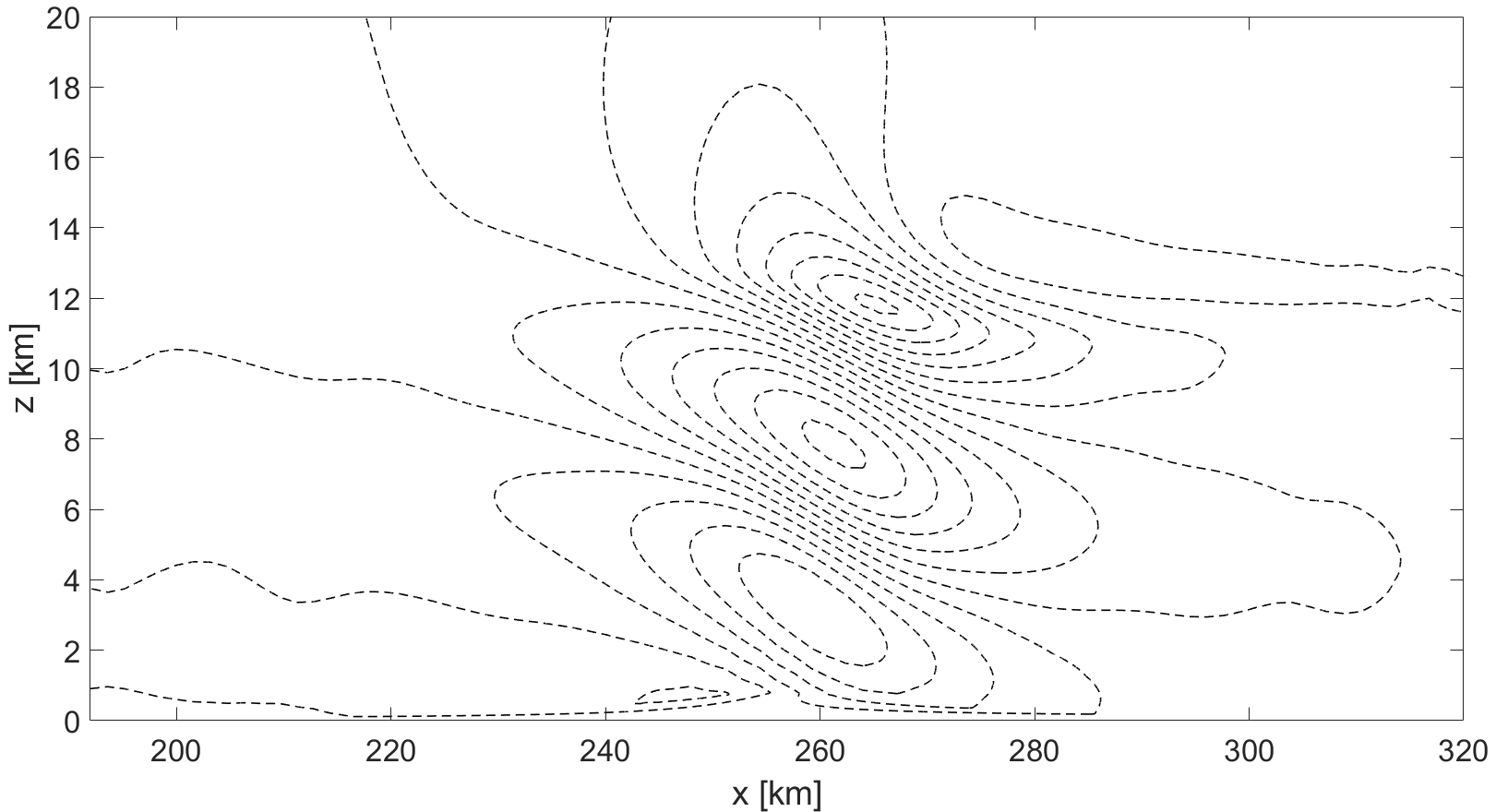} b)
	\caption{Nonlinear hydrostatic flow over a hill at \(t = T_{f}\), a) horizontal velocity deviation, values between \(\SI{-23}{\meter\per\second}\) and \(\SI{28}{\meter\per\second}\) with contour interval of \(\SI{2}{\meter\per\second}\), b) vertical velocity, values between \(\SI{-3.9}{\meter\per\second}\) and \(\SI{3.5}{\meter\per\second}\) with contour interval of \(\SI{0.5}{\meter\per\second}\).}
	\label{fig:nonlinear_hydro_contours}
\end{figure}

\begin{figure}[h!]
	\centering
	\includegraphics[width=0.7\textwidth]{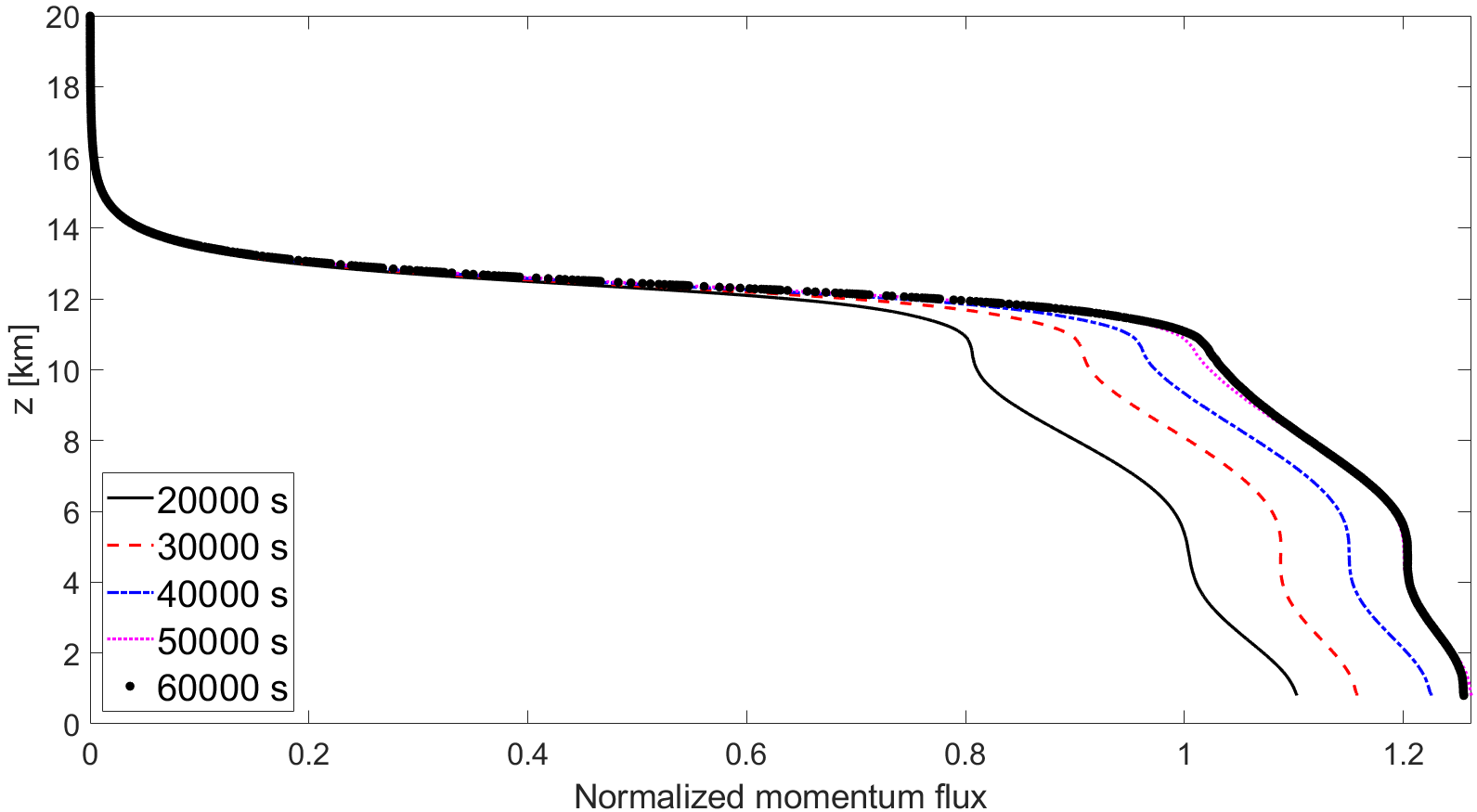}
	\caption{Nonlinear hydrostatic flow over a hill, normalized momentum flux evolution}
	\label{fig:nonlinear_hydro_momentum}
\end{figure}

We consider now the more challenging, nonlinear hydrostatic case considered in \cite{bonaventura:2000, pinty:1995}. The computational domain is \(\Omega = \left(0, 512\right) \times \left(0, 28\right) \hspace{0.05cm} \SI{}{\kilo\meter}\) with \(h_{m} = \SI{800}{\meter}, x_{c} = \SI{256}{\kilo\meter}\) and \(a_{c} = \SI{16}{\kilo\meter}\). The final time is \(T_{f} = \SI{60000}{\second}\). The damping layer is applied starting from \(z = \SI{11.5}{\kilo\meter}\) and in the first and last \(\SI{172}{\kilo\meter}\) along the horizontal direction. The background velocity is \(\overline{u} = \SI{32}{\meter\per\second}\) and the Brunt-V{\"{a}}is{\"{a}}l{\"{a}} frequency \(N\) is equal to \(\SI{0.02}{\per\second}\). The mesh is composed by \(160 \times 112\) elements with \(r = 2\), yielding a resolution of \(\SI{1600}{\meter}\) along \(x\) and \(\SI{125}{\meter}\) along $z, $ whereas the time step is equal to \(\SI{10}{\second}\), yielding a maximum Courant number \(C \approx 1.41\) and \(C_{u} \approx 0.25\). Figure \ref{fig:nonlinear_hydro_contours} shows the contour plots of both the horizontal velocity perturbation and vertical velocity, which compare well with those  presented e.g. in \cite{pinty:1995}. The behaviour over time of the normalized momentum flux is reported in Figure \ref{fig:nonlinear_hydro_momentum} and its value at the surface at \(t = T_{f}\) is approximately equal to \(1.22\), which is comparable to the one obtained in \cite{pinty:1995}. The momentum flux differs from the analytical one because we are no more in a linear regime. Moreover, as explained in \cite{durran:1983}, it is strongly dependent on the position of the absorbing layer. These results confirm the stability and the accuracy of the proposed numerical scheme also in presence
of orography. 

\begin{figure}[h!]
	\centering
	\includegraphics[width=0.9\textwidth]{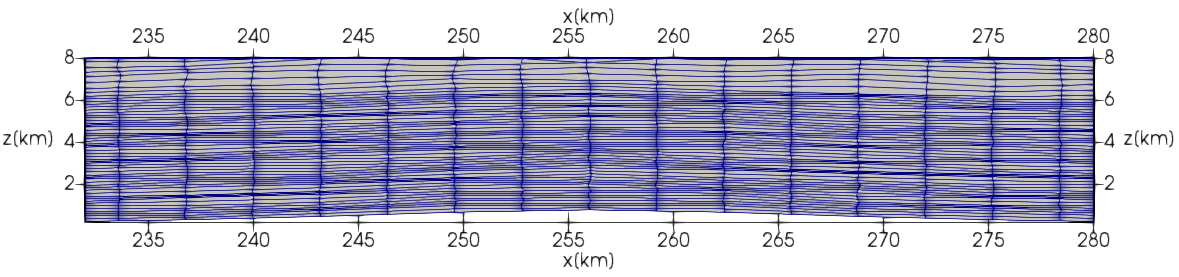}
	\caption{Nonlinear hydrostatic flow over a hill, unstructured grid.}
	\label{fig:nonlinear_hydro_unstructured_mesh}
\end{figure}

\begin{figure}[h!]
	\centering
	\includegraphics[width=0.45\textwidth]{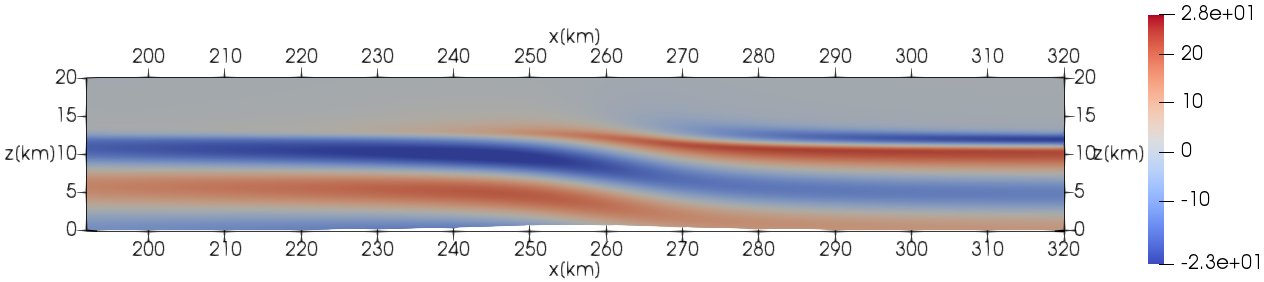} a)
	\includegraphics[width=0.45\textwidth]{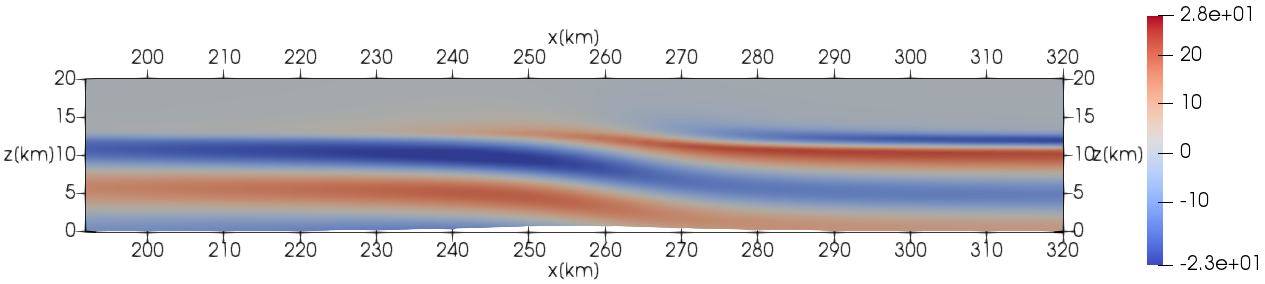} b)
	\caption{Nonlinear hydrostatic flow over a hill at \(t = \SI{50000}{\second}\), a) horizontal velocity deviation with the uniform grid, b) horizontal vertical deviation with the unstructured grid.}
	\label{fig:nonlinear_hydro_unstructured}
\end{figure}

The formulation and the implementation proposed are applicable also to general meshes. We consider now for this benchmark an unstructured grid of quadrilaterals with smaller diameter close to the soil. The mesh is reported in Figure \ref{fig:nonlinear_hydro_unstructured_mesh} and is composed by 5440 elements. Figure \ref{fig:nonlinear_hydro_unstructured} shows also a comparison of the contour plots of the horizontal velocity deviation between the structured grid and the unstructured grid. One can easily notice that a good agreement between the two simulations is established.

\subsection{2D nonhydrostatic flow over a hill}
\label{ssec:nonhydrostatic}

\begin{figure}[h!]
	\centering
	\includegraphics[width=0.7\textwidth]{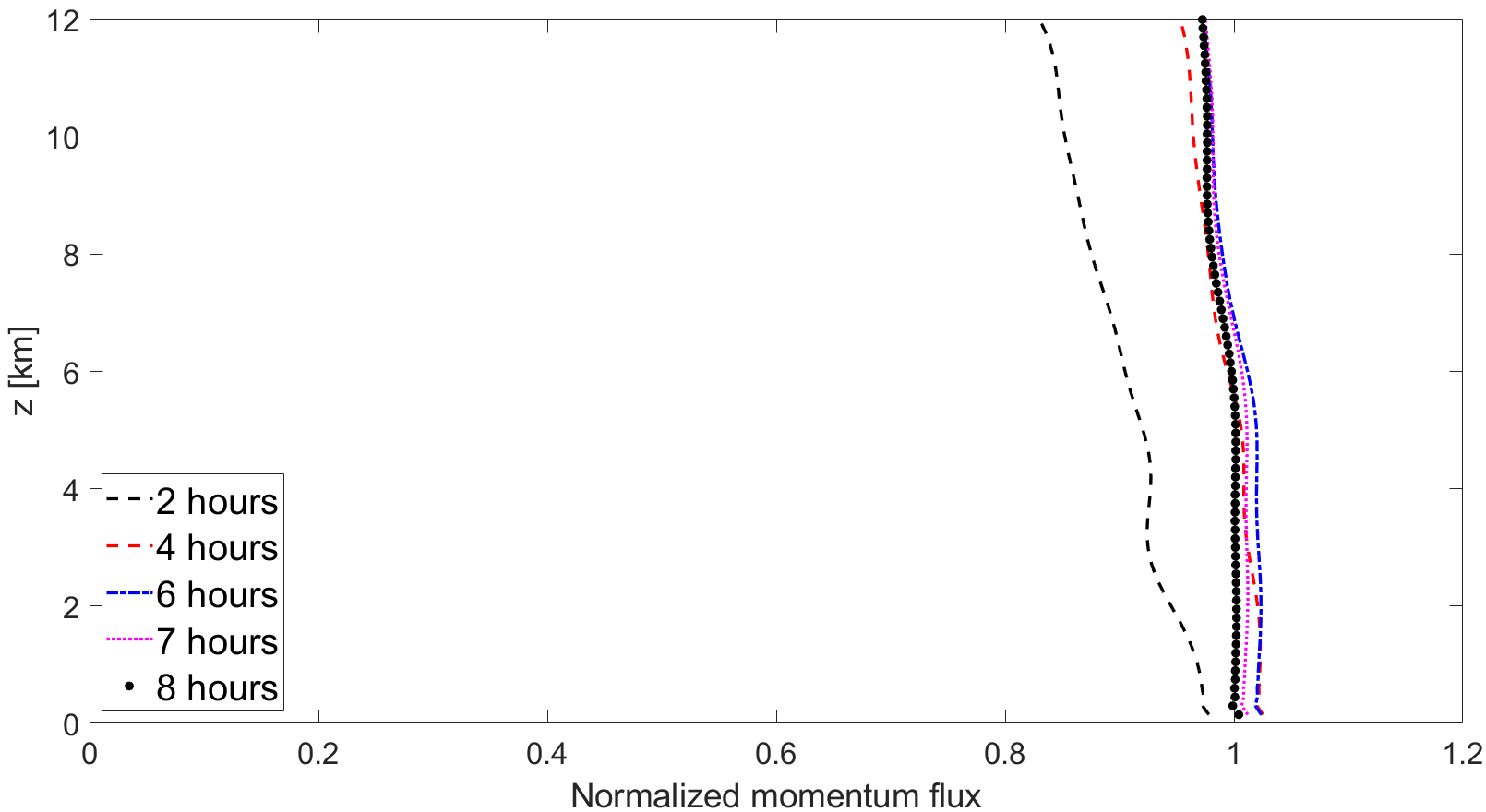}
	\caption{Linear nonhydrostatic flow over a hill, evolution of normalized momentum flux.}
	\label{fig:linear_nonhydro_momentum}
\end{figure}

In this Section, we consider the nonhydrostatic regime, characterized by \(\frac{N a_{c}}{\overline{u}} \approx 1\). The bottom boundary is again described by the function \eqref{eq:versiera_Agnesi}. We first adopt the linear nonhydrostatic configuration described e.g. in \cite{giraldo:2008}. The computational domain is \(\Omega = \left(0, 144\right) \times \left(0, 30\right) \hspace{0.05cm} \SI{}{\kilo\meter}\) with \(h_{m} = \SI{1}{\meter}, x_{c} = \SI{72}{\kilo\meter}\) and \(a_{c} = \SI{1}{\kilo\meter}\). The final time is \(T_{f} = \SI{28800}{\second}\). The initial state of the atmosphere is described by the following potential temperature and Exner pressure, respectively:
$$
\overline{\theta} = \theta_{ref}\exp\left(\frac{N^2}{g}z\right) \ \ \ \ \ \ \
\overline{\pi}    = 1 + \frac{g^2}{c_{p}\theta_{ref}N^2}\left[\exp\left(-\frac{N^2}{g}z\right) - 1\right],
$$
with \(\theta_{ref} = \SI{280}{\kelvin}\) and \(N = \SI{0.01}{\per\second}\). The background velocity \(\overline{u}\) is equal to \(\SI{10}{\meter\per\second}\). Following \cite{klemp:1983}, the analytical momentum flux is given by
$ m^{NH} = 0.457 m^{H} $
and this value will be used to compute the normalized momentum flux for the present case. Wall boundary conditions are applied on the bottom boundary and non-reflecting boundary conditions are employed on the top and lateral boundaries with \(\overline{\lambda}\) such that \(\overline{\lambda}\Delta t = 0.15\). The damping layer is applied in the topmost \(\SI{14}{\kilo\meter}\) of the domain and in the first and last \(\SI{40}{\kilo\meter}\) along the horizontal direction. The mesh is composed by \(200 \times 50\) elements with \(r = 4\), yielding a resolution of \(\SI{180}{\meter}\) along \(x\) and \(\SI{150}{\meter}\) along z, whereas the time step is equal to \(\SI{1}{\second}\), leading to \(C \approx 2.02\) and \(C_{u} \approx 0.06\). Figure \ref{fig:linear_nonhydro_momentum} reports the time evolution of the normalized momentum flux and, as for the linear hydrostatic case in Section \ref{ssec:hydrostatic}, the analytical value is approached as the simulation reaches the steady state.

\begin{figure}[h!]
\centering
	\includegraphics[width=0.45\textwidth]{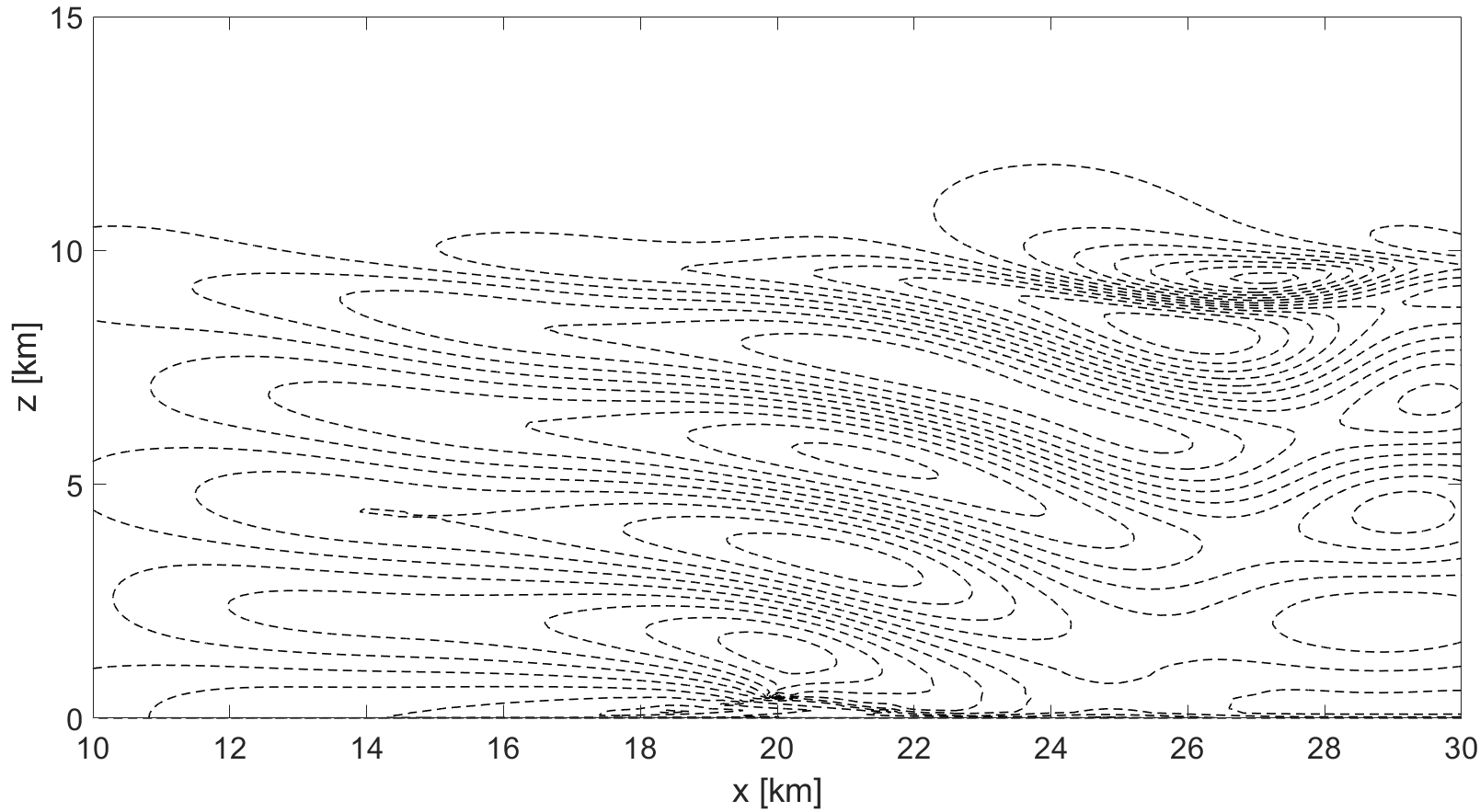} a)
	\includegraphics[width=0.45\textwidth]{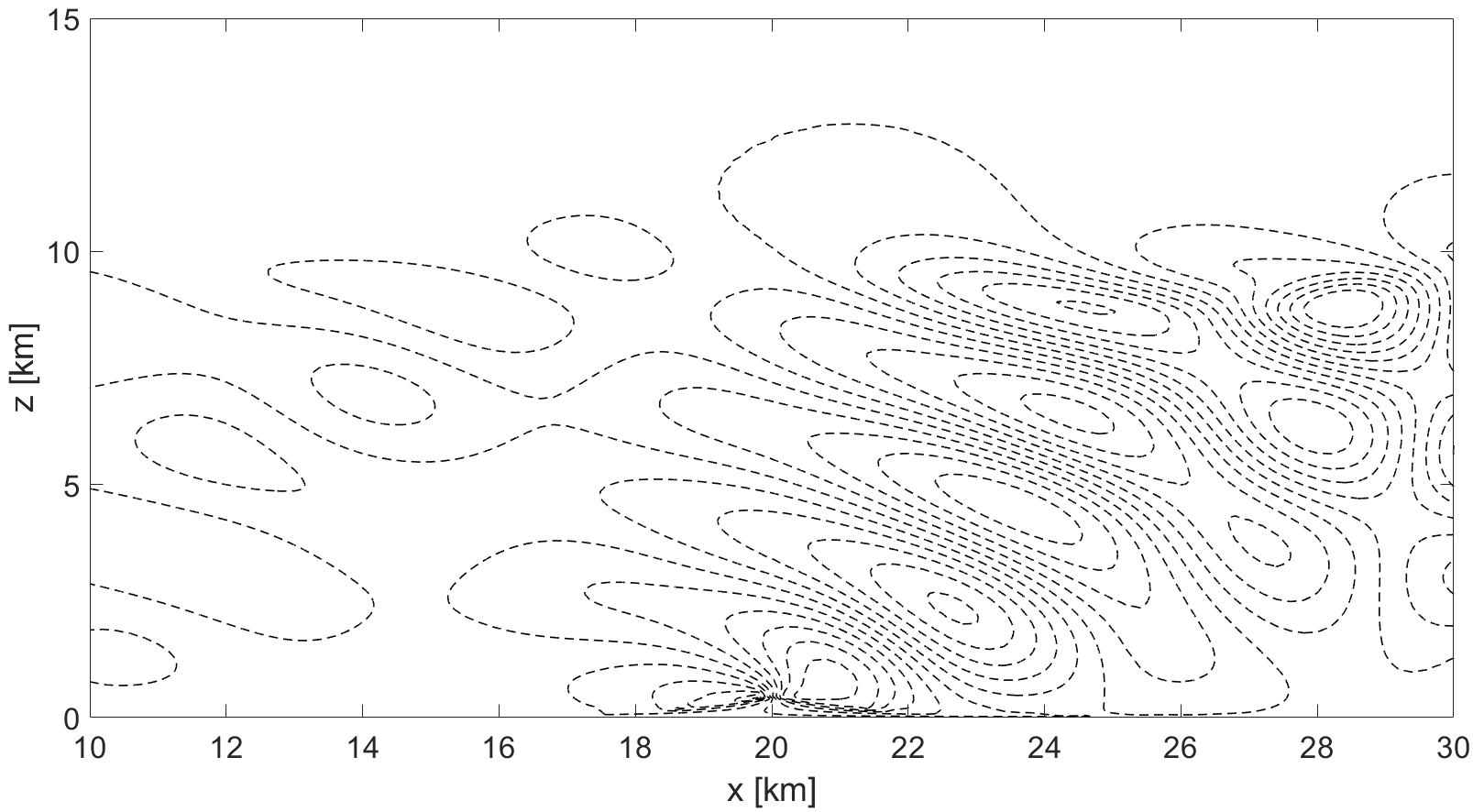} b)
	\includegraphics[width=0.7\textwidth]{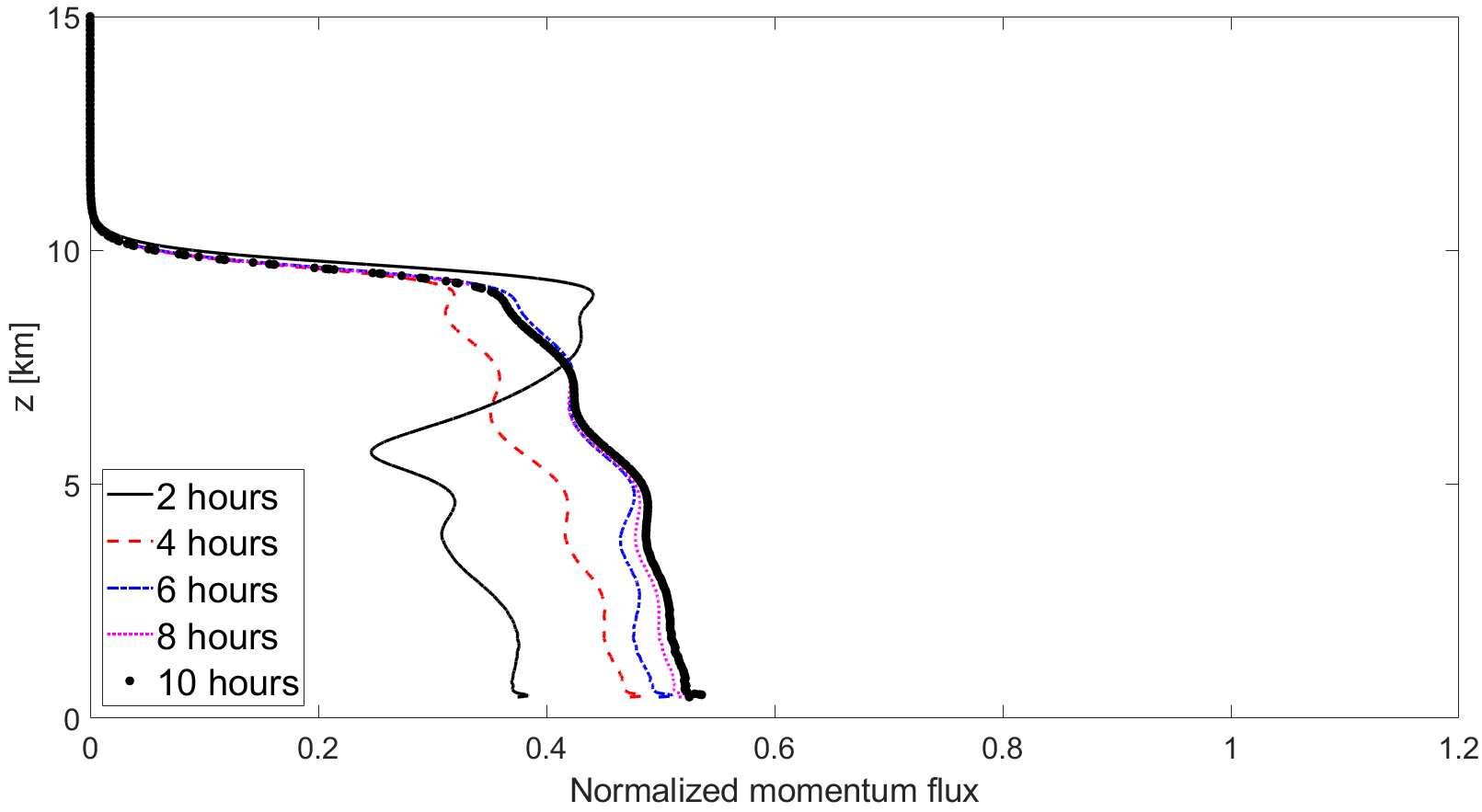} c)
	\caption{Nonlinear nonhydrostatic flow over a hill, a) horizontal velocity deviation, values between \(\SI{-7.2}{\meter\per\second}\) and \(\SI{9.0}{\meter\per\second}\) with contour interval of \(\SI{1.16}{\meter\per\second}\) at \(t = T_{f}\), b) vertical velocity, values between \(\SI{-4.2}{\meter\per\second}\) and \(\SI{4.0}{\meter\per\second}\) with contour interval of \(\SI{0.586}{\meter\per\second}\) at \(t = \SI{18000}{\second}\), c) normalized momentum flux evolution.} 
	\label{fig:nonlinear_nonhydro}
\end{figure}

We consider a nonlinear nonhydrostatic case, see e.g. \cite{tumolo:2015}. The computational domain is \(\Omega = \left(0, 40\right) \times \left(0, 20\right) \hspace{0.05cm} \SI{}{\kilo\meter}\) with \(h_{m} = \SI{450}{\meter}, x_{c} = \SI{20}{\kilo\meter}, a_{c} = \SI{1}{\kilo\meter}, T_{f} = \SI{36000}{\second}, N = \SI{0.02}{\per\second}, \theta_{ref} = \SI{273}{\kelvin}\) and \(\overline{u} = \SI{13.28}{\meter\per\second}\). The damping layer is applied in the topmost \(\SI{11}{\kilo\meter}\) of the domain and in the first and last \(\SI{10}{\kilo\meter}\) along the horizontal direction. The mesh is composed by \(50 \times 50\) elements with \(r = 4\), yielding a resolution of \(\SI{200}{\meter}\) along \(x\) and \(\SI{100}{\meter}\) along $z. $ The time step is equal to \(\SI{0.5}{\second}\), leading to a maximum Courant number \(C \approx 1.13\) and \(C_{u} \approx 0.08\). Figure \ref{fig:nonlinear_nonhydro} shows the contour plots of both horizontal velocity perturbation and vertical velocity, which are analogous to those reported in \cite{tumolo:2015}, as well as the time evolution of the normalized vertical momentum flux. Notice that the momentum flux is normalized by the analytical value \eqref{eq:analytic_momentum_hydro}.

\begin{figure}[h!]
	\centering
	\includegraphics[width=0.48\textwidth]{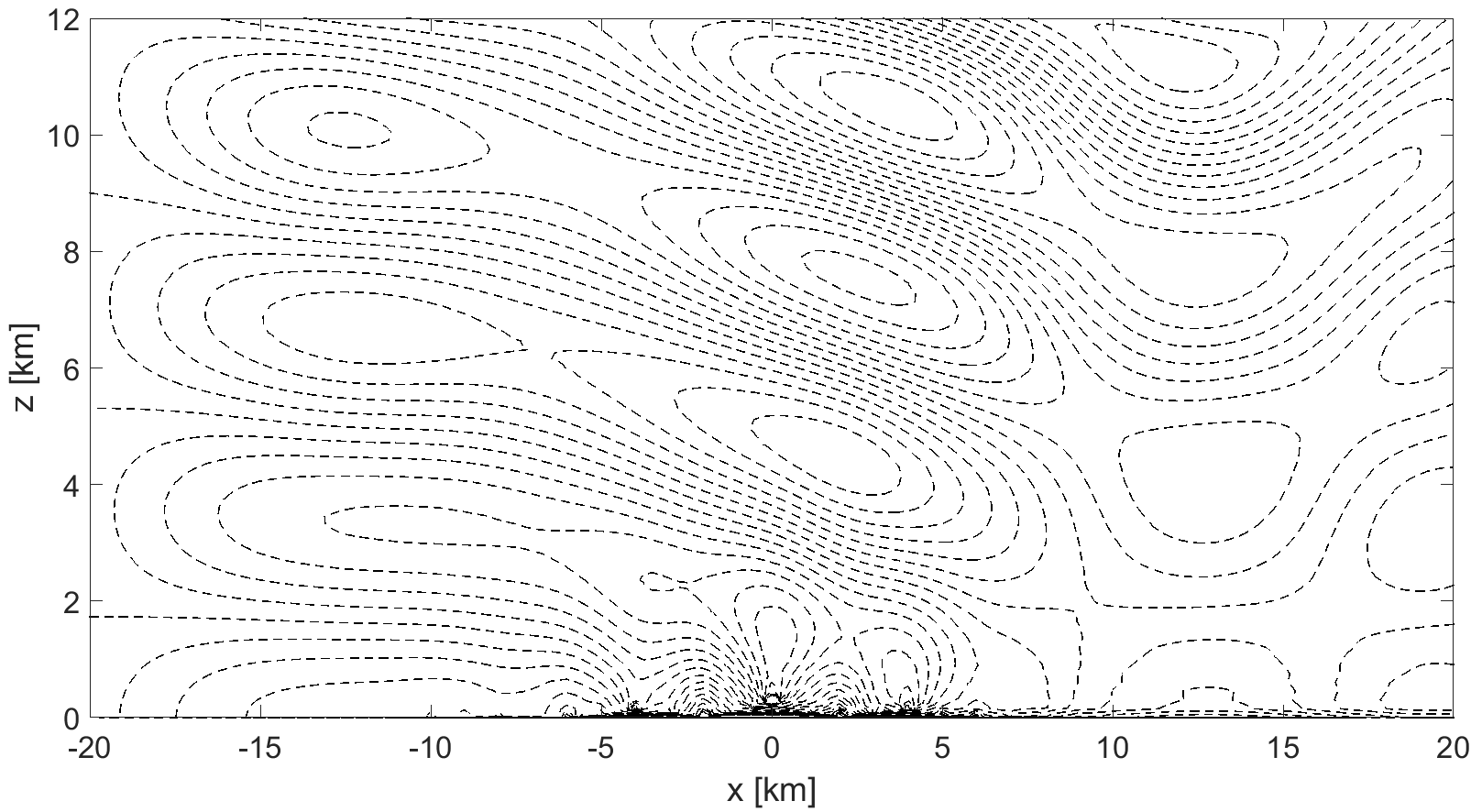} a)
	\includegraphics[width=0.48\textwidth]{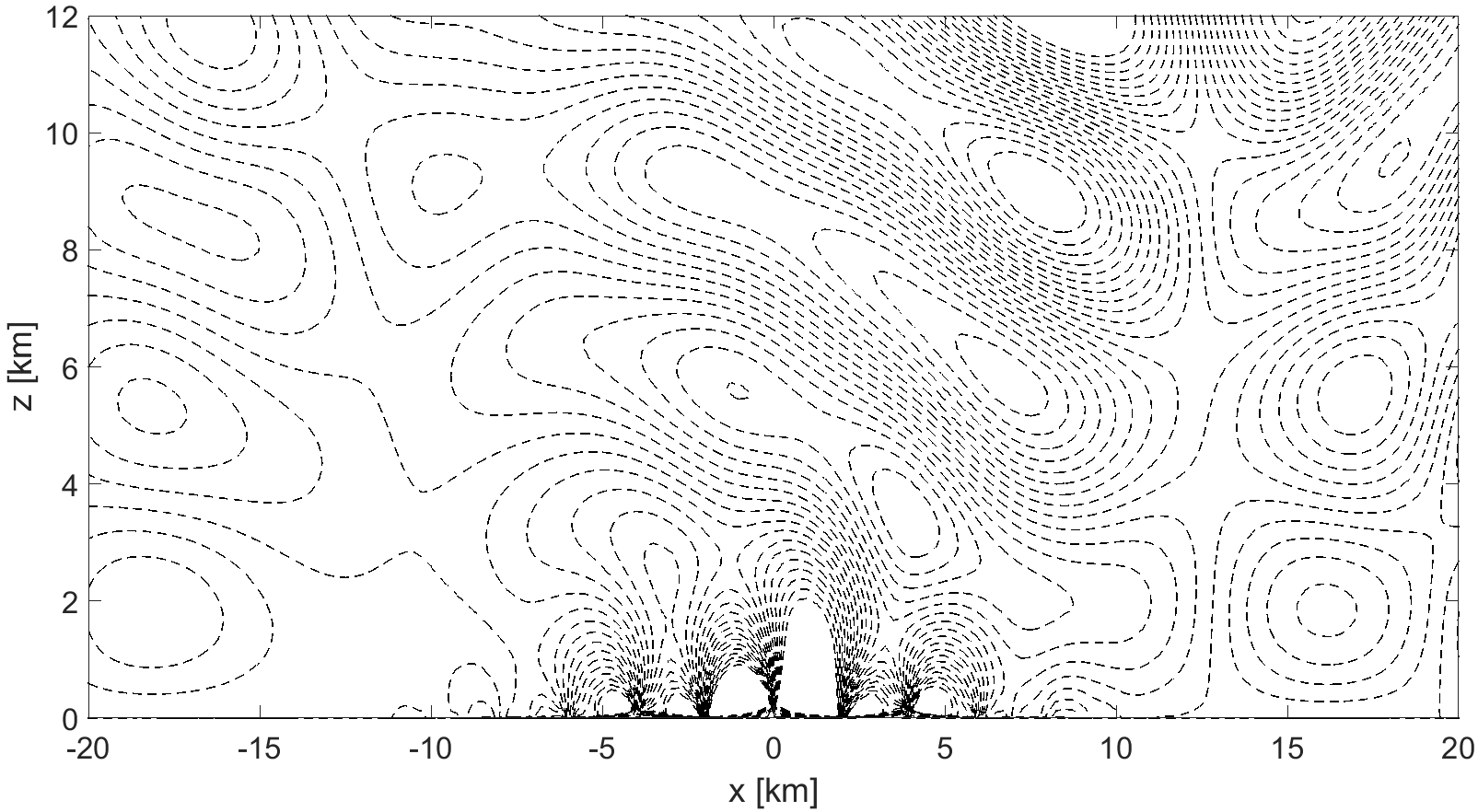} b)
	\caption{Sch{\"{a}}r mountain test case at \(t = T_{f}\), a) horizontal velocity deviation. Contour values are between \(\SI{-2}{\meter\per\second}\) and \(\SI{2}{\meter\per\second}\) with an interval equal to \(\SI{0.2}{\meter\per\second}\), b) vertical velocity. Contour values are between \(\SI{-0.5}{\meter\per\second}\) and \(\SI{0.5}{\meter\per\second}\) with an interval equal to \(5 \cdot 10^{-2} \hspace{0.05cm}\SI{}{\meter\per\second}\).}
	\label{fig:schar_contours}
\end{figure}

As a further nonhydrostatic test, we consider the well known Sch{\"{a}}r mountain test, which consists of a steady-state flow over a five-peak mountain chain \cite{melvin:2019, schar:2002}. The domain is \(\Omega = \left(-50, 50\right) \times \left(0, 30\right) \hspace{0.05cm} \SI{}{\kilo\meter}\) with surface temperature \(T_{ref} = \SI{288}{\kelvin}\), constant buoyancy frequency \(N = \SI{0.01}{\per\second}\) and a background wind \(\bar{u} = \SI{10}{\meter\per\second}\). The mountain profile is defined as
\begin{equation}
h(x) = h_{m}\exp\left[\left(-\frac{x}{a_{c}}\right)^2\right]\cos^2\left(\frac{\pi x}{\lambda_{c}}\right),
\end{equation}
with \(h_m = \SI{250}{\meter}, a_{c} = \SI{5}{\kilo\meter}\) and \(\lambda_{c} = \SI{4}{\kilo\meter}\). The background density and pressure have the same expression as in Section \ref{ssec:igw}, with \(\theta_{ref} = \SI{288}{\kelvin}\) and the final time is \(T_f = \SI{18000}{\second}\). The damping layer is applied in the topmost \(\SI{10}{\kilo\meter}\) of the domain and in the first and last \(\SI{10}{\kilo\meter}\) along the horizontal direction with \(\overline{\lambda}\Delta t = 1.2\). The mesh is composed by \(100 \times 50\) elements with \(r = 4\), leading to a resolution of \(\SI{250}{\meter}\) along \(x\) and of \(\SI{150}{\meter}\) along \(z\), whereas the time step is equal to \(\SI{2.5}{\second}\), yielding a maximum acoustic Courant number \(C \approx 2.02\) and a maximum advective Courant number \(C_{u} \approx 0.09\). Figure \ref{fig:schar_contours} shows the contour plots of both horizontal velocity perturbation and vertical velocity, which are analogous to those reported in \cite{giraldo:2008, melvin:2019}. 

\begin{figure}[h!]
	\centering
	\includegraphics[width=0.48\textwidth]{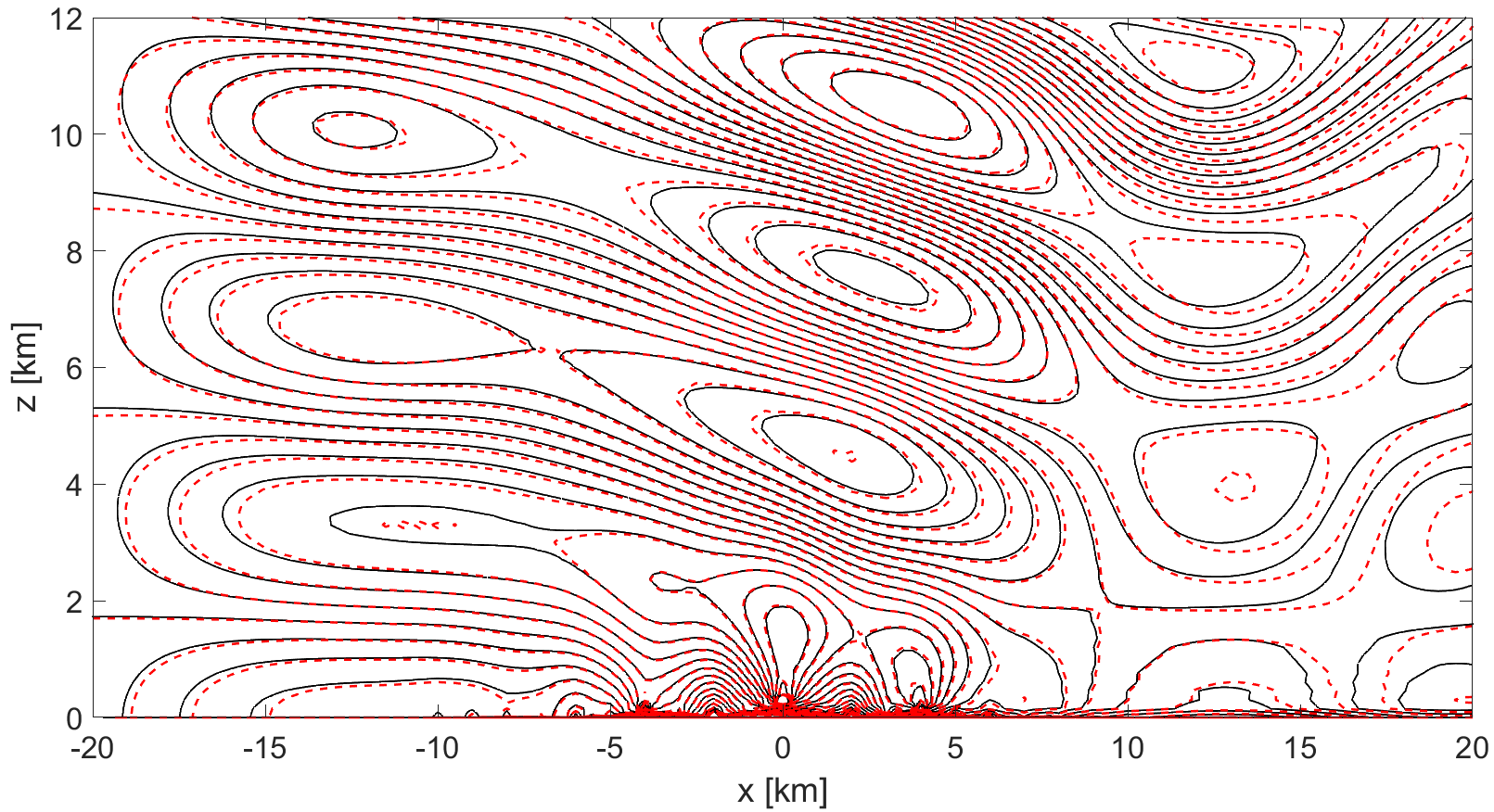} a)
	\includegraphics[width=0.48\textwidth]{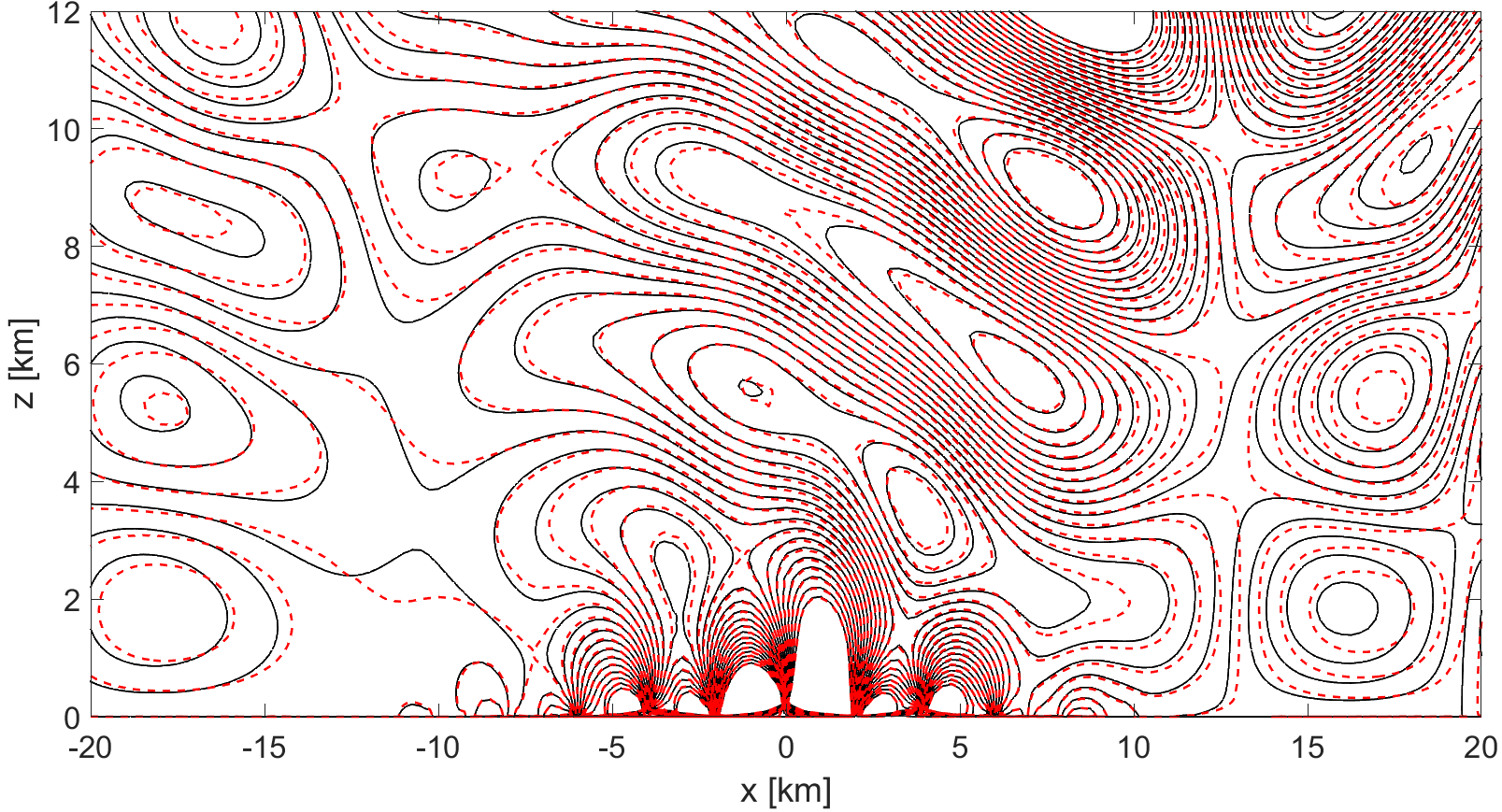} b)
	\caption{Sch{\"{a}}r mountain test case at \(t = T_{f}\), a) horizontal velocity deviation. Contour values are between \(\SI{-2}{\meter\per\second}\) and \(\SI{2}{\meter\per\second}\) with an interval equal to \(\SI{0.2}{\meter\per\second}\), b) vertical velocity. Contour values are between \(\SI{-0.5}{\meter\per\second}\) and \(\SI{0.5}{\meter\per\second}\) with an interval equal to \(5 \cdot 10^{-2} \hspace{0.05cm}\SI{}{\meter\per\second}\). The black lines denote the solution with uniform grid, whereas the red lines denote the solution with the non-conforming grid.}
	\label{fig:schar_contours_adaptive}
\end{figure}

We consider now for this test case a non-conforming mesh refinement over the soil, applying a coarsening above \(z = \SI{2}{\kilo\meter}\), so as to obtain a resolution of \(\SI{500}{\meter}\) along the horizontal direction and \(\SI{300}{\meter}\) along the vertical direction. The grid is composed by \(1550\) elements. Figure \ref{fig:schar_contours_adaptive} shows a comparison of both horizontal and the vertical velocity deviation contours between the results obtained using the non-conforming grid and those obtained with the uniform mesh. A computational saving time of around \(25 \%\) is achieved by the use of the non-conforming grid and one can easily notice a good agreement of the two solutions, without spurious reflections at the interface between the two grids. This opens the way to the use of non-conforming grids for flows over orography. We plan in future work to fully exploit the potentialities of this technique and to show that larger scale features can be adequately resolved on coarser meshes and an increased resolution can be used only locally over the orography without affecting the correct	momentum transfer in the vertical direction.

\subsection{3D medium-steep bell-shaped hill}
\label{ssec:3D_mountain}

\begin{figure}[h!]
	\centering
	\includegraphics[width=0.5\textwidth]{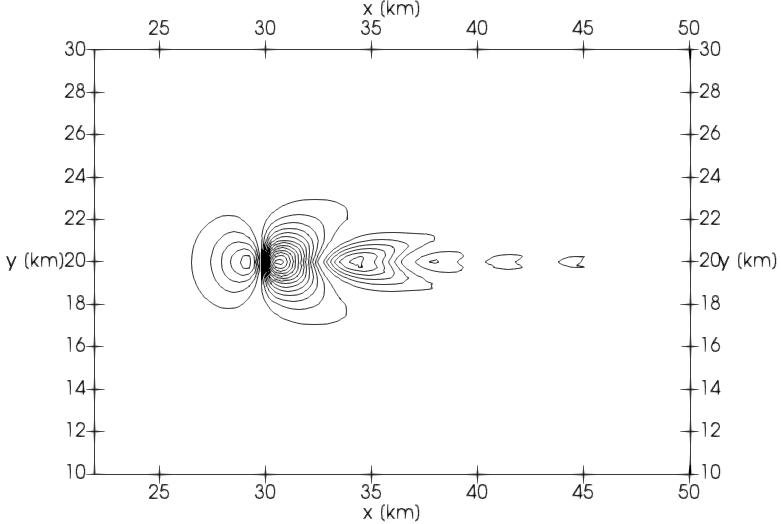}
	\caption{3D mountain benchmark, \(x - y\) slice at \(z = \SI{800}{\meter}\) of the vertical velocity. The values are between \(\SI{-1.5}{\meter\per\second}\) and \(\SI{1.3}{\meter\per\second}\) with contour interval of \(\SI{0.1}{\meter\per\second}\)}
	\label{fig:3D_hill}
\end{figure}

In this Section, we consider the three-dimensional flow over a bell-shaped hill, see e.g. \cite{melvin:2019}. The computation domain is \(\Omega = \left(0, 60\right) \times \left(0, 40\right) \times \left(0, 20\right) \hspace{0.05cm} \SI{}{\kilo\meter}\). The mountain profile is defined as
\begin{equation}\label{eq:versiera_Agnesi_3D}
h(x, y) = \frac{h_{m}}{\left[1 + \left(\frac{x - x_{c}}{a_{c}}\right)^2 + \left(\frac{y - y_{c}}{a_{c}}\right)^2\right]^{3/2}},
\end{equation}
with \(h_m = \SI{400}{\meter}, a_{c} = \SI{1}{\kilo\meter}, x_{c} = \SI{30}{\kilo\meter}\) and \(y_{c} = \SI{20}{\kilo\meter}\). We consider as buoyancy frequency \(N = \SI{0.01}{\per\second}\) and a background velocity \(\overline{u} = \SI{10}{\meter\per\second}\). We are therefore in a nonhydrostatic regime since \(\frac{N a_{c}}{\overline{u}} = 1\). The background density and pressure have the same expression of as in  Section \ref{ssec:igw}, with \(\theta_{ref} = \SI{293.15}{\kelvin}\) and the final time is \(T_f = \SI{3600}{\second}\). The damping layer is applied in the topmost \(\SI{6}{\kilo\meter}\) of the domain and in the first and last \(\SI{20}{\kilo\meter}\) along the lateral boundaries with \(\overline{\lambda}\Delta t = 1.2\). In order to increase the resolution around the hill, we consider a non-uniform grid by taking a resolution of \(\SI{250}{\meter}\) between \(x = \SI{25}{\kilo\meter}\) and \(x = \SI{40}{\kilo\meter}\) and a resolution of \(\SI{250}{\meter}\) between \(y = \SI{12}{\kilo\meter}\) and \(y = \SI{28}{\kilo\meter}\). A uniform resolution of \(\SI{500}{\meter}\) is considered along the vertical direction \(z\), as well as for the remaining part of the lateral boundaries. The mesh is composed by 8288 elements with polynomial degree \(r = 4\), whereas the time step is equal to \(\SI{2}{\second}\), yielding a maximum acoustic Courant number \(C \approx 1.95\) and a maximum advective Courant number \(C_{u} \approx 0.1\). The results in Figure \ref{fig:3D_hill} are in agreement with those reported in \cite{melvin:2019}. 

\section{Conclusions}
\label{sec:conclu}

We have presented an application of the IMEX discretization of the compressible Euler equations proposed in \cite{orlando:2022b} to numerical modelling of atmospheric flows. The method includes non-conforming $h-$ refinement, as implemented in the framework provided by the numerical library \textit{deal.II} \cite{arndt:2022, bangerth:2007}. The use of an open source library enhances code portability and embeds the development of codes specific for NWP applications in  a software environment that is being continuously adapted to novel architectures and extended to include novel versions of the DG method. Simulations of classical benchmarks demonstrate that the method can simulate accurately small scale flows in presence of gravity and idealized flows over orography. In future developments, we aim to show the potential of non-conforming mesh refinement to increase locally the resolution over complex orography. Moreover, we will consider the inclusion of more complex physical phenomena, such as for example turbulence, water vapour transport and adiabatic heating, in order to demonstrate that all the typical features of a high resolution numerical weather prediction model can be included in the proposed adaptive framework without significant loss of accuracy or efficiency.

\section*{Acknowledgements}
We would like to thank all the reviewers, who have provided several constructive comments that helped to improve the presentation of the topics discussed in the paper. The simulations have been partly run at CINECA thanks to the computational resources made available through the ISCRA-C projects  SIDICoNS - HP10CLPLXI  and   NUMNETF - HP10C06Y02. This work has been partly supported by the ESCAPE-2 project, European Union’s Horizon 2020 Research and Innovation Programme (Grant Agreement No. 800897).

\bibliographystyle{cas-model2-names}
\bibliography{DG_IMEX_ATMO}

\end{document}